\documentclass[a4paper,12pt]{article}
\usepackage{booktabs}
\usepackage[toc,page]{appendix}

\pdfoutput=1

\usepackage[paper=letterpaper,margin=1in]{geometry}

\usepackage{amsmath,amssymb,amsfonts,epsfig,cite,setspace,bigstrut,longtable,array,breqn,url,color}
\usepackage{cleveref}
\usepackage{caption,setspace}
\usepackage{tabularx}

\usepackage{pdfpages,lipsum}

\usepackage[OT2,T1]{fontenc}
\DeclareSymbolFont{cyrletters}{OT2}{wncyr}{m}{n}
\DeclareMathSymbol{\Sha}{\mathalpha}{cyrletters}{"58}

\begin{document}

%%%%%%%%%%%%%%%%%%%%%%%%%%%%%%%%%%%%%%%%%%%%%%%%%%%%%%%%%%%%%%%%%%

%\includepdf[pages={12-15,23,45-49}]{main.pdf}

\vskip 0.25in

\newcommand{\sref}[1]{\S~\ref{#1}}
\newcommand{\nn}{\nonumber}
\newcommand{\tr}{\mathop{\rm Tr}}
\newcommand{\ord}{\mathop{\rm Ord}}
\newcommand{\slog}{\mathop{\rm sLog}}
\newcommand{\sgn}{\mathop{\rm sgn}}

\newcommand{\comment}[1]{}

\newcommand{\cM}{{\cal M}}
\newcommand{\cW}{{\cal W}}
\newcommand{\cN}{{\cal N}}
\newcommand{\cH}{{\cal H}}
\newcommand{\cK}{{\cal K}}
\newcommand{\cY}{{\cal Y}}
\newcommand{\cZ}{{\cal Z}}
\newcommand{\cO}{{\cal O}}
\newcommand{\cA}{{\cal A}}
\newcommand{\cB}{{\cal B}}
\newcommand{\cC}{{\cal C}}
\newcommand{\cD}{{\cal D}}
\newcommand{\cE}{{\cal E}}
\newcommand{\cF}{{\cal F}}
\newcommand{\cX}{{\cal X}}
\newcommand{\IA}{\mathbb{A}}
\newcommand{\IP}{\mathbb{P}}
\newcommand{\IQ}{\mathbb{Q}}
\newcommand{\IH}{\mathbb{H}}
\newcommand{\IR}{\mathbb{R}}
\newcommand{\IC}{\mathbb{C}}
\newcommand{\IF}{\mathbb{F}}
\newcommand{\IV}{\mathbb{V}}
\newcommand{\II}{\mathbb{I}}
\newcommand{\IZ}{\mathbb{Z}}
\newcommand{\re}{{\rm~Re}}
\newcommand{\im}{{\rm~Im}}

\newcommand{\tmat}[1]{{\tiny \left(\begin{matrix} #1 \end{matrix}\right)}}
\newcommand{\mat}[1]{\left(\begin{matrix} #1 \end{matrix}\right)}

% make bibliography single-spaced
\let\oldthebibliography=\thebibliography
\let\endoldthebibliography=\endthebibliography
\renewenvironment{thebibliography}[1]{%
\begin{oldthebibliography}{#1}%
\setlength{\parskip}{0ex}%
\setlength{\itemsep}{0ex}%
}%
{%
\end{oldthebibliography}%
}

\newtheorem{theorem}{\bf THEOREM}
\def\thetheorem{\thesection.\arabic{theorem}}
\newtheorem{proposition}{\bf PROPOSITION}
\def\thetheorem{\thesection.\arabic{proposition}}
\newtheorem{observation}{\bf OBSERVATION}
\def\thetheorem{\thesection.\arabic{observation}}
\newtheorem{conjecture}{\bf CONJECTURE}
\def\thetheorem{\thesection.\arabic{conjecture}}

\def\theequation{\thesection.\arabic{equation}}
\newcommand{\setall}{\setcounter{equation}{0}
        \setcounter{theorem}{0}}
\newcommand{\setequation}{\setcounter{equation}{0}}
\renewcommand{\thefootnote}{\fnsymbol{footnote}}

~\\
\vskip 1cm

\begin{center}
{\Large \bf  %Machine Learning meets Pure Maths: \\ 
Machine Learning meets Number Theory: \\ The Data Science of Birch-Swinnerton-Dyer}

\medskip
\vspace{4mm}

{\large Laura Alessandretti$^{1,2}$, Andrea Baronchelli$^{2,3}$, Yang-Hui He$^{2,4,5}$}%\blfootnote{}

\vspace{1mm}

\renewcommand{\arraystretch}{0.5} 
{\small
{\it
\begin{tabular}{rl}
  ${}^{1}$ & Copenhagen Center for Social Data Science,\\
 &University of Copenhagen, Copenhagen K, 1353, Denmark \\
  ${}^{2}$ &
  Department of Mathematics, City, University of London, EC1V 0HB, UK\\
  ${}^{3}$ & The Alan Turing Institute, London NW1 2DB, UK \\
  ${}^{4}$ & Merton College, University of Oxford, OX14JD, UK\\
  ${}^{5}$ &
  School of Physics, NanKai University, Tianjin, 300071, P.R.~China
  
\end{tabular}
}
~\\
~\\
~\\
l.alessandretti@gmail.com, \
Andrea.Baronchelli.1@city.ac.uk, \ 
hey@maths.ox.ac.uk\\
}
\renewcommand{\arraystretch}{1.5} 

\end{center}

\vspace{5mm}

\begin{abstract}
Empirical analysis is often the first step towards the birth of a conjecture. This is the case of the Birch-Swinnerton-Dyer (BSD) Conjecture describing the rational points on an elliptic curve, one of the most celebrated unsolved problems in mathematics. Here we extend the original empirical approach, to the analysis of the Cremona database of quantities relevant to BSD, inspecting more than 2.5 million elliptic curves by means of the latest techniques in data science, machine-learning and topological data analysis. 

Key quantities such as rank, Weierstrass coefficients, period, conductor, Tamagawa number, regulator and order of the Tate-Shafarevich group give rise to a high-dimensional point-cloud whose statistical properties we investigate. We reveal patterns and distributions in the rank versus Weierstrass coefficients, as well as the Beta distribution of the BSD ratio of the quantities.
Via gradient boosted trees, machine learning is applied in finding inter-correlation amongst the various quantities. We anticipate that our approach will spark further research on the statistical properties of large datasets in Number Theory and more in general in pure Mathematics.

%We apply the latest techniques in data-science, notably machine-learning (ML) and topological data analysis (TDA) to the Cremona database of quantities relevant to the Birch-Swinnerton-Dyer (BSD) Conjecture on around 2.5 million elliptic curves.
%This empirical analysis is very much in the spirit of the original investigations of BSD which gave birth to the conjecture.
%Key quantities such as rank, Weierstrass coefficients, period, conductor, Tamagawa number, regulator and order of the Tate-Shafarevich group give rise to a high-dimensional point-cloud whose statistical properties we investigate.
%We present patterns and distributions in the rank versus Weierstrass coefficients, as well as the Beta distribution of the BSD ratio of the quantities.
%Via gradient boosted trees, ML is applied in finding inter-correlation amongst the various quantities.

\end{abstract}

%\end{titlepage}

\newpage

\tableofcontents

%%%%%%%%%%%%%%%%%%
\section{Introduction and Summary}
Elliptic curves $\cE$ occupy a central stage in modern mathematics, their geometry and arithmetic providing endless insights into the most profound structures.
The celebrated Conjecture of Birch and Swinnerton-Dyer \cite{BSD} is the key result dictating the behaviour of $\cE$ over finite number fields and thereby, arithmetic.
Despite decades of substantial progress, the proof of the conjecture remains elusive.
To gain intuition, a highly explicit and computational programme had been pursued by Cremona \cite{cremona}, in cataloguing all elliptic curves up to isogeny and expressed in canonical form, to conductors into the hundreds of thousands.

Interestingly, a somewhat similar situation exists for the higher dimensional analogue of elliptic curves considered as Ricci-flat K\"ahler manifolds, viz., Calabi-Yau manifolds.
Though Yau \cite{yau} settled the Calabi-Yau Conjecture \cite{calabi}, much remains unknown about the landscape of such manifolds, even over the complex numbers. For instance, even seemingly simple questions of whether there is a finite number of topological types of Calabi-Yau $n$-folds for $n\geq 3$ is not known -- even though it is conjectured so.
Nevertheless, because of the pivotal importance of Calabi-Yau manifolds to superstring theory, theoretical physicists have been constructing ever-expanding datasets thereof over the last few decades (cf.~\cite{He:2018jtw} for a pedagogical introduction).

Given the recent successes in the science of ``big data'' and machine learning, it is natural to examine the database of Cremona \cite{lmfdb} using the latest techniques of Data Science. 
Indeed, such a perspective has been undertaken for Calabi-Yau manifolds and the landscape of compactifications in superstring theory in high-energy physics, ranging from machine-learning \cite{He:2017aed} to statistics \cite{Douglas:2003um,He:2015fif}. Indeed, \cite{He:2017aed,Krefl:2017yox,Carifio:2017bov,Ruehle:2017mzq} brought about a host of new activities in machine-learning within string theory; moreover, \cite{He:2017aed,He:2018jtw} and the subsequent work in \cite{Bull:2018uow,Altman:2018zlc,He:2019vsj,Jejjala:2019kio,He:2019nzx,Brodie:2019dfx,Ashmore:2019wzb} introduced the curious possibility that machine-learning should be applied to at least stochastically avoid expensive algorithms in geometry and combinatorics and to raise new conjectures.

Can artificial intelligence help with understanding the syntax and semantics of mathematics?
While such profound questions are better left to the insights of Turing and Voevodsky, the more humble question of using machine-learning to recognizing patterns which might have been missed by standard methods should be addressed more immediately.
Preliminary experiments such as being able to ``guess''- to over 99\% accuracy and confidence - the ranks of cohomology groups without exact-sequence-chasing (having seen tens of thousands of examples of known bundle cohomologies) \cite{He:2018jtw}, or whether a finite group is simple without recourse to theorem of Noether and Sylow (having been trained on thousands of Cayley tables) \cite{He:2019nzx} already point to this potentiality.

In our present case of number theory, extreme care should, of course, be taken.
Patterns in the primes can be notoriously deceptive, as exemplified by the likes of Skewe's constant.
Indeed, a sanity check to let neural networks predict the next prime number in \cite{He:2017aed} yielded a reassuring null result.
Nevertheless, one should not summarily disregard all experimentation in number theory: after all, the best neural network of the 19th century - the mind of Gau\ss\ - was able to pattern-spot $\pi(x)$ to raise the profound Prime Number Theorem years before the discovery of complex analysis to allow its proof.

The purpose of this paper is to open up dialogue between data scientists and number theorists, as the aforementioned works have done for the machine-learning community with geometers and physicists, using Cremona's elliptic curve database as a concrete arena.
The organization is as follows.
We begin with a rapid introduction in Section 2, bearing in mind of the diversity of readership, to elliptic curves in light of BSD. In Section 3, we summarize Cremona's database and perform preliminary statistical analyses beyond simple frequency count. Subsequently, Section 4 is devoted to machine-learning various aspects of the BSD quantities and Section 5, to their topological data analyses and persistent homology.
Finally, we conclude in Section 5 with the key results and discuss prospects for the future.

%%%%%%%%%%%%%%%%%%%
\section{Elliptic Curves and BSD}
Our starting point is the Weierstra\ss\ model
\begin{equation}\label{weier}
y^2 + a_1xy + a_3y = x^3 + a_2 x^2 + a_4 x + a_6
\end{equation}
of an elliptic curve $\cE$ over $\IQ$, where $(x,y) \in \IQ$ and the coefficients $a_i \in \IZ$.
The discriminant $\Delta$ and J-invariant of $\cE$ are obtained in a standard way as
\begin{equation}\label{deltaj}
\Delta(\cE)=-b_2^2b_8 - 8 b_4^3 -27 b_6 ^2 + 9 b_2 b_4 b_6
\ , 
\qquad
j(\cE)= \frac{c_4^3}{\Delta}
\end{equation}
where
$b_2 := a_1^2 + 4 a_2$,
$b_4 := 2a_4 + a_1 a_3$, 
$b_ 6 := a_3^2 + 4 a_6$,
$b_8 : = a_1^2 a_6 + 4 a_2 a_6 - a_1 a_3 a_4 + a_2 a_3^2 - a_4^2$ and
$c_4 := b_2^2 - 24b_4$.
Smoothness is ensured by the non-vanishing of $\Delta$ and isomorphism (isogeny) between two elliptic curves, by the equality of $j$.

An algorithm of Tate and Laska \cite{tate,laska}
\footnote{
In particular, consider the transformation between the coefficiens $a_i$ and $a_i'$ between two elliptic curves $\cE$ and $\cE'$:
\begin{align*}
        u  a_1' &= a_1+2s,                               \\
        u^2a_2' &= a_2-sa_1+3r-s^2,                      \\
        u^3a_3' &= a_3+ra_1+2t,                          \\
        u^4a_4' &= a_4-sa_3+2ra_2-(t+rs)a_1+3r^2-2st,    \\
        u^6a_6' &= a_6+ra_4+r^2a_2+r^3-ta_3-t^2-rta_1.   \\
\end{align*}
for $u,s,t \in \IQ$, then relating the points $(x,y) \in \cE$ and $(x',y') \in \cE'$ as
\[
 x = u^2x' + r          \ , \quad
 y = u^3y' + su^2x' + t \\        
\]
yields $u^{12}\Delta' = \Delta$ and hence $j' = j$, and thus the isomorphism.
}
can then be used to bring the first 3 coefficients $a_{1,2,3}$ to be $0, \pm 1$, transforming \eqref{weier} into a {\em minimal Weierstra\ss\ model}.
%https://mathoverflow.net/questions/10607/writing-down-minimal-Weierstra\ss\ -equations
%minimal also in minimizing Delta
Thus, for our purposes, an elliptic curve $\cE$ is specified by the pair of integers $(a_4, a_6)$ together with a triple $(a_1, a_2, a_3) \in \{-1, 0, 1\}$.
From the vast subject of elliptic curves, we will need this 5-tuple, together with arithmetic data to be presented in the ensuing.

\subsection{Rudiments on the Arithmetic of $\cE$}
%We have a neither the space nor the expertise to give an account of the rich subject of the theory of elliptic curves.
This subsection serves as essentially a lexicon for the quantities which we require, presented in brief, itemized form.
\begin{description}

\item[Conductor and Good/Bad Reduction: ]
The conductor is the product over all (finite many) number of primes $p$ - the primes of bad reduction -
where $\cE$ reduced modulo $p$ becomes singular (where $\Delta = 0$).
All other primes are called good reduction.

\item[Rank and Torsion: ]
The set of rational points on $\cE$ has the structure of an Abelian group, $\cE(\IQ) \simeq \IZ^r \times T$.
The non-negative integer $r$ is called the rank, its non-vanishing would signify an infinite number of rational points on $\cE$.
The group $T$ is called the torsion group and can be only one of 15 groups by Mazur's celebrated theorem \cite{mazur}, viz.,
the cyclic group $C_n$ for $1 \leq n \leq 10$ and $n=12$, as well as the direct product $C_2 \times C_n$ for $n=2,4,6,8$. 

\item[L-Function and Conductor: ]
The Hasse-Weil Zeta-function of $\cE$ can be defined, given a finite field $\IF_{q = p^n}$, as the generating functions
$z_p$ (the local) and $Z$ (the global):
\begin{align}
\nn
Z_p(t; \cE) & := \exp\left(  \sum\limits_{n=1}^\infty \frac{ \cE (\IF_{p^n})}{n} t^n \right) \ , \\
\label{zeta}
Z(s; \cE) & := \prod\limits_p Z_p( t := p^{-s}; \cE) \ .
\end{align}
Here, in the local zeta function $Z_p(t; \cE)$, $\cE (\IF_{p^n})$ is the number of points of $\cE$ over the finite field  
and the product is taken over all primes $p$  to give the global zeta function $Z(s)$.

The definition \eqref{zeta} is applicable to general varieties, and for elliptic curves, the global zeta function simplifies (cf.~\cite{silverman}) to a product of the Riemann zeta function $\zeta$ and a so-called L-function as
\begin{equation}
Z(s; \cE) = \frac{\zeta(s)\zeta(s-1)}{L(s;\cE)} \ ,
\end{equation}
where
\begin{align}
\nn
L(s; \cE) &= \prod\limits_p L_p(s; \cE)^{-1} \ , \quad
L_p(s; \cE) := \left\{
\begin{array}{lcl}
(1-\alpha_{p}p^{-s}+p^{1-2s}), &&p\nmid N\\
(1-\alpha_{p}p^{-s}),&& p\mid N{\text{ and }}p^{2}\nmid N\\
1,&&p^{2}\mid N
\end{array}
\right. \ .
\end{align}
In the above, $\alpha_p = p+1 -$ counts the number of points of  $\cE$ mod $p$ for primes of good reduction and $\pm1$ depending on type of bad reduction.
The positive integer $N$ which controls, via its factorization, these primes, is the conductor of $\cE$.
 
Importantly, the L-function has analytic continuation  \cite{TW} to $\IC$ so that the variable $s$ is not a merely dummy variable like $t$ in the local generating function, but renders $L(s; \cE)$ a complex analytic function.

\item[Real Period: ]
The periods of a complex variety is usually defined to be the integral of some globally defined holomorphic differential over a basis in homology.
Here, we are interested in the real period, defined as (using the minimal Weierstra\ss\  model)
\begin{equation}
\IR \ni \Omega := \int_{\cE(\IR)} |\omega| \ , \quad \omega = \frac{dx}{2y + a_1 x + a_3}
\end{equation}
over the set of real points $\cE(\IR)$ of the elliptic curve.

\item[Tamagawa Number: ]
Now, $\cE$ over any field is a group, thus in particular we can define $\cE(\IQ_p)$ over the $p$-adic field $\IQ_p$ for a given prime, as well as its subgroup $\cE^0(\IQ_p)$ of points which have good reduction.
We define the index in the sense of groups
\begin{equation}
c_p := \Big[\cE(\IQ_p) : \cE(\IQ_p)\Big] \ ,
\end{equation}
which is clearly equal to 1 for primes of good reduction since then $\cE^0(\IQ_p) = \cE(\IQ_p)$.
The Tamagawa number is defined to be the product over all primes of bad-reduction of $c_p$, i.e.,
\begin{equation}
\text{Tamagawa Number } = \prod\limits_{p \mid N} c_p \ .
\end{equation}

\item[Canonical Height: ]
For a rational point $P = \frac{a}{b}$, written in minimal fraction form with $\gcd(a,b) = 1$, $a,b \in \IZ$ and $b >0$, we can define a naive height
$h(P) := \log \max(|a|, b)$. Then, a {\em canonical height} can be defined as 
\begin{equation}
\hat{h}(P) = \lim\limits_{n \to \infty} n^{-2} h(n P) \ ,
\end{equation}
where $n P = P + \ldots + P$ ($n$-times) is the addition of under the group law of $\cE$.
This limit exists and renders $\hat{h}$ the unique quadratic form on $\cE(\IQ) \otimes \IR$ such that $\hat{h} - h$ is bounded.
An explicit expression of $\hat{h}(P)$ in terms of $a,b$ can be found, e.g., in Eq4 and Eq5 of \cite{BGZ}.
The canonical height defines a natural bilinear form
\begin{equation}
2 \left<P,P'\right> = \hat{h}(P + P') - \hat{h}(P) -\hat{h}(P')
\end{equation}
for two points $P,P' \in \cE(\IQ)$ and as always, $P+P'$ is done via the group law.

\item[Regulartor: ]
Given the infinite (free Abelian) part of $\cE(\IQ)$, viz., $\IZ^r$, let its generators be  $P_1, \ldots, P_r$, then we can define the regulator
\begin{equation}
R_{\cE} = \det \left<P_i, P_j\right> \ , \quad  i,j = 1,\ldots, r
\end{equation}
where the pairing is with the canonical height and defines an $r \times r$ integer matrix.
For $r=0$, $R$ is taken to be 1 by convention.

\item[Tate-Shafarevich Group: ]
Finally, one defines group cohomologies $H^1(\IQ, \cE)$ and $H^1(\IQ_p, \cE)$ between which there is a homomorphism
(cf.~e.g., \cite{milne}{IV.2} for a detailed description).
We can then define the Tate-Shafarevich Group $\Sha$ of $\cE$ as the kernel of the homomorphism
\begin{equation}
\Sha(\cE) := \ker \left(
	H^1(\IQ, \cE) \longrightarrow \prod\limits_p  H^1(\IQ_p, \cE)
\right) \ .
\end{equation}
This is the most mysterious part of the arithmetic of elliptic curves, it is conjectured to be a finite Abelian group.
For ranks $r=0,1$, this has been proven (cf.~the survey of \cite{rs}) but in general this is not known.

\end{description}

\subsection{The Conjecture}\label{s:conj}
With the above definitions, we can at least present the celebrated
\begin{conjecture}[Birch--Swinnerton-Dyer (Weak Version)]
The order of the zero of $L(s; \cE)$ at $s = 1$ is equal to the rank $r$, 
\[
\ord\limits_{s \to 1} L(s; \cE) = r(\cE) \ .
\]
That is, the Taylor series around 1 is $L(s; \cE) \sim c (s-1)^r$ with some complex coefficient $c$.
\end{conjecture}
In fact, a stronger version of the conjecture predicts precisely what the Taylor coefficient $c$ should be:
\begin{conjecture}[Birch--Swinnerton-Dyer (Strong Version)]
The Taylor coefficient of $L(s; \cE)$ at $s=1$ is given in terms of the regulator $R$, Tamagawa number $\prod\limits_{p \mid N} c_p$, (analytic) order of the Tate-Shafarevich group $\Sha$, 
and the order of the torsion group $T$.
Specifically, let $ L(s; \cE) = \sum\limits_{r} \frac{L^{(r)}(1; \cE)}{r!} (s-1)^r$, then
\[
 \frac{L^{(r)}(1; \cE)}{r!} =
 \frac{|\Sha| \cdot \Omega \cdot R \cdot \prod\limits_{p \mid N} c_p}{|T|^2}
 \ .
\]
\end{conjecture}

%%%%%%%%%
\section{Elliptic Curve Data}\setall
BSD arose from extensive computer experimentation, the earliest of its kind, by Birch and Swinnerton-Dyer.
Continuing along this vein, Cremona \cite{cremona} then compiled an impressive of list of 2,483,649 isomorphism classes of elliptic curves over $\IQ$ and explicitly computed the relevant quantities introduced above.
This is available freely online at \cite{lmfdb}.
\subsection{Cremona Database}
The database of Cremona,  on which the rest of the paper will focus, associates to each minimal Weierstra\ss\ model (given by the coefficients $(a_1,a_2,a_3) \in \{-1,0,1\}$ and $(a_4,a_6) \in \IZ$; generically, these last two coefficients have very large magnitude), the following:
\begin{itemize}
\item the conductor $N$, ranging from 1 to 400,000;
\item the rank $r$, ranging from 0 to 4;
\item the torsion group $T$, whose size ranges from 1 to 16;
\item the real period $\Omega$, a real number ranging from approximately $2.5 \cdot 10^{-4}$ to $6.53$.
\item the Tamagawa number $\prod\limits_{p \mid N} c_p$, ranging from 1 to 87040;
\item the order of the Tate-Shafarevich group (exactly when known, otherwise given numerically), ranging from 1 to 2500;
\item the regulator $R \in \IZ_{>0}$, ranging from approximately $0.01$ to $3905.84$.
\end{itemize}

A typical entry, corresponding to the curve
$y^2 + x y = x^{3} -  x^{2} - 453981 x + 117847851$ (labelled as ``314226b1'' and with $\Delta = 2 \cdot 3^{3} \cdot 11 \cdot 23^{8}$ and $j=2^{-1} \cdot 3^{3} \cdot 11^{-1} \cdot 23 \cdot 199^{3}$ which are readily computed from\eqref{deltaj}) would be
\begin{align}
\nn
(a_1,a_2,a_3,a_4,a_6) & = (1,-1,0,- 453981,117847851) \qquad \Longrightarrow\\
& \left\{
\begin{array}{rcl}
N &=& 314226 = 2 \cdot 3^3 \cdot 11 \cdot 23^2 \\
r &=& 0 \\
R &=& 1 \\
\Omega &\simeq& 0.56262 \\
\prod\limits_{p \mid N} c_p &=& 3 \\
|T| &=& 3 \\
|\Sha| &=& 1 \ . \\
\end{array}
\right.
\end{align}

%%%%
\subsection{Weierstra\ss\ Coefficients}
Let us begin with a statistical analysis of the minimal Weierstra\ss\ coefficients themselves.
It so happens that in the entire database, there are only 12 different sets of values of $(a_1,a_2,a_3)$, we tally all of the curves in the following histogram, against rank and $(a_1,a_2,a_3)$:
\[
\begin{array}{|c|c|c|c|c|c|}\hline
 \text{$(a_1,a_2,a_3) \backslash $Rank} & 0 & 1 & 2 & 3 & 4 \\ \hline \hline
 \{0,-1,0\} & 126135 & 155604 & 30236 & 659 & 0 \\
 \{0,-1,1\} & 17238 & 24593 & 7582 & 399 & 0 \\
 \{0,0,0\} & 172238 & 213780 & 40731 & 698 & 0 \\
 \{0,0,1\} & 28440 & 39235 & 11187 & 506 & 0 \\
 \{0,1,0\} & 118942 & 157003 & 34585 & 722 & 0 \\
 \{0,1,1\} & 18016 & 27360 & 9609 & 426 & 0 \\
 \{1,-1,0\} & 102769 & 127198 & 25793 & 551 & 1 \\
 \{1,-1,1\} & 96995 & 128957 & 28940 & 604 & 0 \\
 \{1,0,0\} & 66411 & 98092 & 25286 & 612 & 0 \\
 \{1,0,1\} & 71309 & 94595 & 20907 & 548 & 0 \\
 \{1,1,0\} & 69759 & 88403 & 18293 & 496 & 0 \\
 \{1,1,1\} & 67834 & 91717 & 21197 & 458 & 0 \\ \hline
\end{array}
\]
We see that most of the curves are of smaller rank, with only a single instance of $r=4$.
This is in line with the recent result of \cite{BS} that most elliptic curves are rank 1; in fact, over 2/3 of elliptic curves obey the BSD conjecture \cite{BSZ}.

To give an idea of the size of the a-coefficients involved, the largest one involved in the database is
\begin{equation}
\vec{a} =  \{
{1, 0, 0, -40101356069987968, -3090912440687373254444800} 
\} \ ,
\end{equation}
which is of rank 0.

Even though  it is a conjecture that the rank $r$ can be arbitrarily large,
the largest confirmed rank \cite{elkies} so far known in the literature is 19, corresponding (the last term being a 72-digit integer!) to 
\begin{align}
\nn
&a_1 = 1, \
a_2 = 1, \
a_3 = -1, \
a_4 = 31368015812338065133318565292206590792820353345,
\\
\label{rk19}
&a_6 = 
302038802698566087335643188429543498624522041683874493555186062568159847 \ .
\end{align}
This extraordinary result is clearly not in the database due to the large rank.

\begin{figure}[!h!t!b]
%\begin{picture}(500,350)(0,-160)
%\put(60,0){$\chi  = 2(h^{1,1} - h^{2,1})$}
%\put(80,170){$h^{1,1} + h^{2,1}$}
%\put(360,0){$\chi$}
%\put(220,120){\#}
\includegraphics[trim=0mm 0mm 0mm 0mm, clip, width=6in]{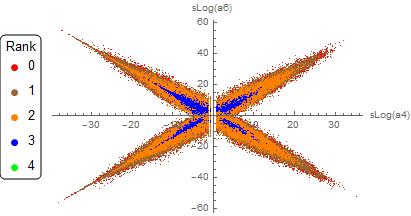}
%\end{picture}
\caption{
{\sf {\small 
A scatter plot of $(\slog(a_4), \ \slog(a_6))$ for all 2,483,649  elliptic curves in the Cremona database.
Different ranks are marked with different colours.
}}
\label{f:sloga4a6}
}
\end{figure}

One of the first steps in data visualization is a {\em principle component analysis} where features of the largest variations are extracted.
The minimal Weierstra\ss\ model gives a natural way of doing so, since the $(a_1, a_2, a_3)$ coefficients take only 12 values and we can readily see scatter plot of $(a_4,a_6)$.
Now, due to the large variation in these coefficients, we define a signed natural logarithm for $x \in \IR$ as
\begin{equation}
\slog(x) = \left\{
\begin{array}{rcl}
\sgn(x) \log (x) \ , && x \neq 0 \\
0 \ , && x=0 \ .
\end{array}
\right.
\end{equation}
We present this scatter plot of $(\slog(a_4), \ \slog(a_6))$ in Fig.~\ref{f:sloga4a6}.
Therein, we plot all the data points (i.e., for all different values of $(a_1, a_2, a_3)$) together, distinguishing rank by colour (rank 4 has only a single point as seen from the table above).

%%%
%%%

The first thing one would notice is the approximate cross-like symmetry, even within rank.
However, this is not trivial because the transformation
\begin{equation}
(a_4, a_6) \longrightarrow (\pm a_4, \pm a_6)
\end{equation}
is by no means a rank preserving map.
For instance, a single change in sign in $a_4$, could result in rank change from 1 to 3:
\begin{equation}
r(\{0, 1, 1, -10, 20\}) = 3 \ ,  \quad r(\{0, 1, 1, +10, 20\}) = 1 \ .
\end{equation}
Examples of a similar nature abound.
The next feature to notice is that the size of the cross shrinks as the rank increases.
This is rather curious since the largest rank case of \eqref{rk19} has far {\em larger} coefficients.
This symmetry is somewhat reminiscent of mirror symmetry for Calabi-Yau 3-folds, where every compact smooth such manifold with Hodge numbers $(h^{1,1}, h^{2,1})$ has a mirror manifold with these values exchanged.

%%%%%%%%%%%%%%%%%
\subsection{Distributions of $a_4$ and $a_6$}

%%%
\begin{figure}[h!]
\centering
\includegraphics[width=\linewidth]{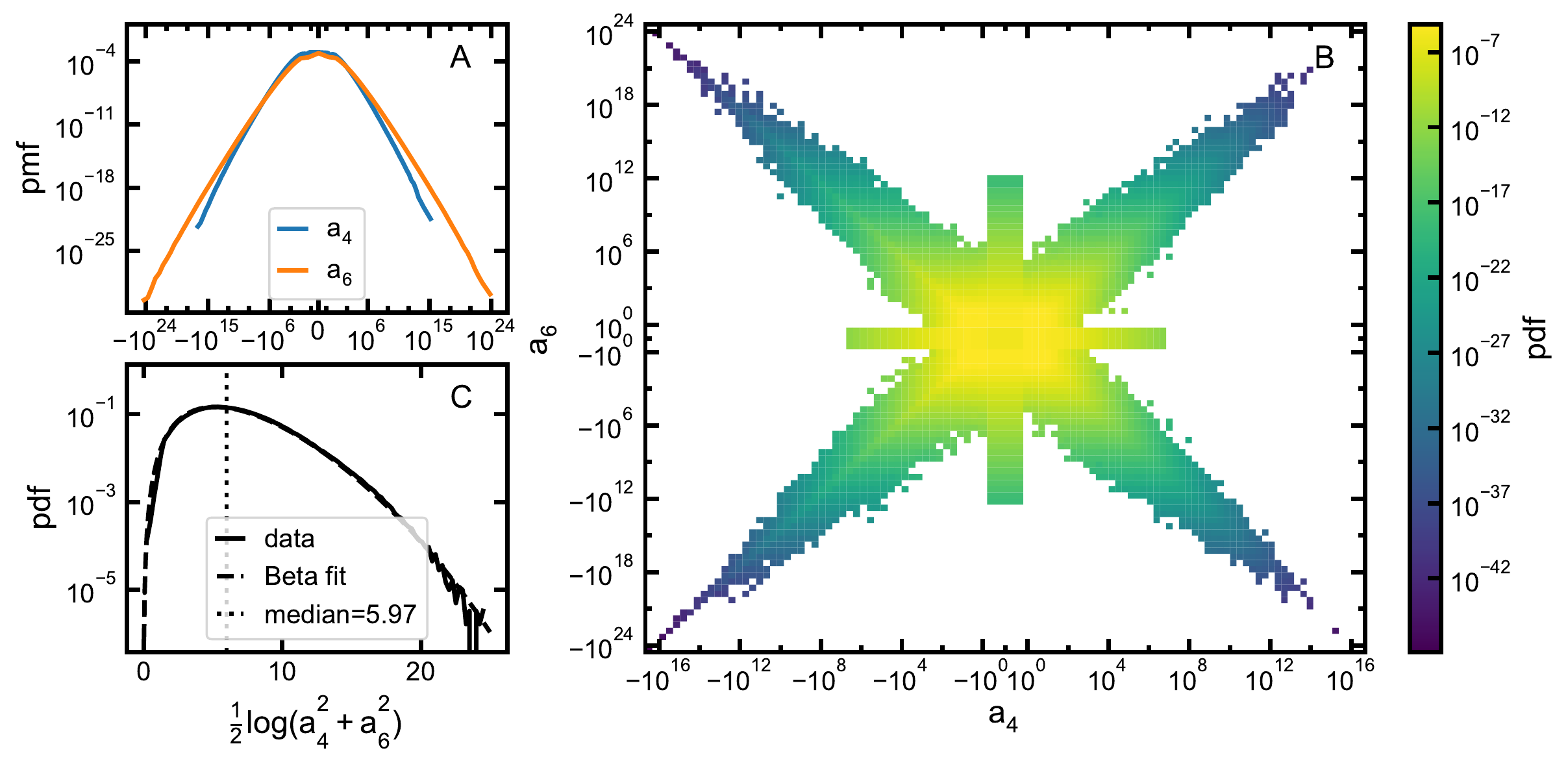}
\caption{\textbf{The distributions of $a_4$ and $a_6$.} (A) Probability mass of $a_4$ (blue line) and $a_6$ (orange line). (B) Joint probability of $a_4$ and $a_6$. The colour indicates the density of points within each bin (see colour-bar). Note that the figure axes are in symlog scale (linear scale between $-1$ and $1$, logarithmic scale for other values in the range). Hence, we can see also the density corresponding to $a_4=0$ and $a_6=0$ (cross in the middle). (C) Probability density distribution for the logarithm(in base 10) of $\sqrt{a_4^2+a_6^2}$ (when $a_4>0$ and $a_6>0$, filled line), the corresponding Beta distribution fit with parameters $\alpha=4.1$, $\beta=25.0$ and $s=44.1$ (dashed line), and the median value (dotted line). }
\label{a4-a6-pdf}
\end{figure}
%%%

Fortified by the above initial observations, let us continued with more refined study of the distribution of $a_4$ and $a_6$.
First, let us plot the distribution of each individually, normalized by the total, i.e., as {\em probability mass functions}.
These are shown in part (A) to the left of Fig~\ref{a4-a6-pdf}. Note that the horizontal axes (binning) is done logarithmically.
We see that the distributions of both $a_4$ and $a_6$ are symmetric with respect to $0$ (Fig.~\ref{a4-a6-pdf}-A), with $a_4$ spanning $\sim 8$ orders of magnitude smaller as compared to $a_6$. This is just to give an idea of the balanced nature of Cremona's data, that elliptic curves with $\pm a_4$ and $\pm a_6$ are all constructed.

Next, in part (B) of Fig~\ref{a4-a6-pdf}, we plot the joint probability mass function of the pair $(a_4,a_6)$ with colour showing the frequency as indicated by the colour-bar to the right. 
We see that, as discussed in Fig.~\ref{f:sloga4a6}, there is a cross-like symmetry. Here, since we are not separating by rank, the symmetry is merely a reflection of the constructions of the dataset, that $\pm a_4$ and $\pm a_6$ are all present.
What is less explicable is that it should be a cross shape and what is the meaning of the boundary curve beyond which there does not seem to be any minimal models. 
For reference, the central rectilinear cross indicates the cases of $a_4=0$ and $a_6=0$ respectively.

Finally, we compute the Euclidean distance $d := \sqrt{a_4^2+a_6^2}$ from the origin and study its probability distribution.
This is shown in part (C) of Fig.~\ref{a4-a6-pdf},
We find that half of the data lies within a radius of $\sim 10^6$ from the origin. The logarithm of $d$ can be well fitted with a Beta probability distribution:
\begin{equation}
f(x, \alpha, \beta, s) = K \cdot \frac{x}{s
}^{\alpha-1} \left(1-\frac{x}{s}\right)^{\beta-1} \ ,
\end{equation}
with parameters $\alpha=4.1$, $\beta=25.0$ and $s=44.1$.
Thus, whilst there are a number of coefficients of enormous magnitude, the majority still have relatively small ones.
%%%
\begin{figure}[t!!!]
\centering
\includegraphics[width=\linewidth]{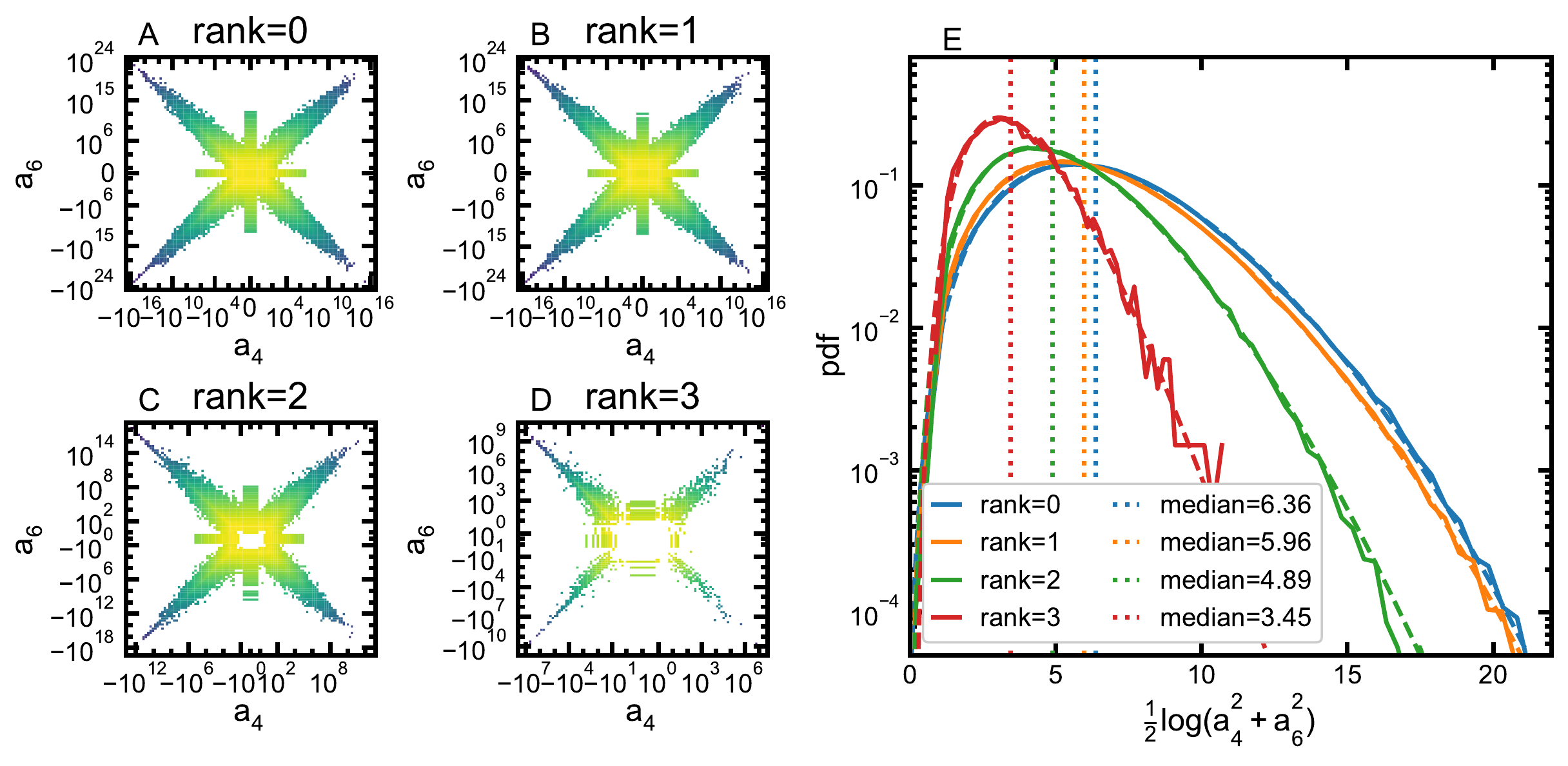}
\caption{\textbf{Different distributions of $a_4$ and $a_6$ for different ranks.} (A-D) Joint probability distribution of $a_4$ and $a_6$ for values of the rank $r=0$ (A), $r=1$ (B), $r=2$ (C), $r=3$ (D). (E) Probability density distribution for the logarithm(in base 10) of $\sqrt{a_4^2+a_6^2}$ (when $a_4>0$ and $a_6>0$, filled lines) for various value of the rank $r$, the corresponding Beta distribution fit (dashed lines), and the corresponding median values (dotted lines).}
\label{a4-a6-rank}
\end{figure}
%%%

%%%%%%%%%%%%%%%%%%%%
\paragraph{Differences by Rank: }
As with our initial observations, we now study the variation of $(a_4,a_6)$ with respect to the rank $r$ of the elliptic curves. 
First, in Fig. \ref{a4-a6-rank} parts A-D, we plot the joint distributions of $a_4$ and $a_6$ for $r=0,1,2,3$ respectively. We can see that they differ signigicantly from each other, under permutation test \cite{permutation} at confidence level $\alpha=0.01$.

Next, Fig. \ref{a4-a6-rank} E show the probability distribution functions for our Euclidean distance $\sqrt{a_4^2+a_6^2}$ for the different ranks.
We find that the median Euclidean distance from the center decreases for higher values of $r$. 
In fact, we see that the median values of $a_4$ and $a_6$ increase with the rank $r$ (see Fig. \ref{boxplots}D-E which we will discuss shortly).
Again, each is individually well-fitted by the Gamma distribution.
In tables \ref{Table1} and  \ref{Table2} in the Appendix, we show some statistics of $a_4$ and $a_6$ including their mean, standard deviation, median, and the number of zero entries, for given rank $r$ and values of $(a_1$, $a_2$, $a_3)$.

%%%%%%%%%%%%%%%%%%%%
\subsection{Distributions of various BSD Quantities}

%%%%
\begin{figure}[h!!!]
\centering
\includegraphics[width=\linewidth]{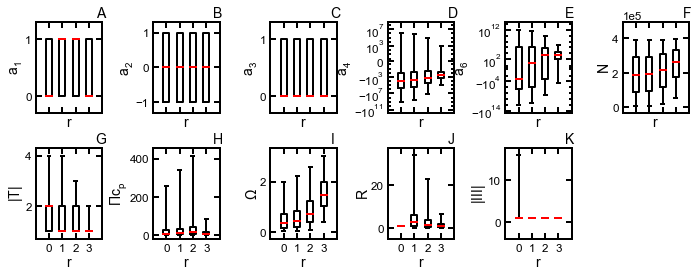}
\caption{\textbf{Characteristics of the elliptic curves based on their rank.} Boxplots of $a_1$, $a_2$, $a_3$, $a_4$, $a_6$ 
in parts (A)-(E) respectively;
boxplots of
$N$, $|T|$, $\prod\limits_{p \mid N} c_p$, $\Omega$, $R$ and $|\Sha|$ in parts (F)-(K) respectively,
all for different values of the rank $r=0,1,2,3$. The red line shows the median value, the boxes enclose  $50\%$ of the distribution, and the whiskers, $95\%$ of the distribution.}
\label{boxplots}
\end{figure}
%%%%%
Now, the coefficients $a_i$ are the inputs of the dataset, each specifying a minimal model of an elliptic curve, and to each should be associated the numerical tuple $(r, N, R, \Omega, \prod\limits_{p \mid N} c_p, |T|, |\Sha|)$ for the rank, the conductor, the regulator, the real period, the Tamagawa number, the order of the torsion group and the order of Tate-Shafarevich group, respectively.
It is now expedient to examine the distribution of ``output'' parametres.
\comment{
%%%%
\subsection{Conductor}
The conductor has a range from 1 to $4 \times 10^5$, while the Tate-Shafarevich group has size which spreads from 
1 to 2500.
In Fig.~\ref{f:NSha}, we plot the sLog of the conductor versus that of the order of $\Sha$..
}

As always, we arrange everything by rank $r=0,1,2,3$ and in Fig.~\ref{boxplots} show the box plots of the variation around the median (drawn in red).
The boxes enclose  $50\%$ of the distribution, and the whiskers, $95\%$ of the distribution.
We see that, as detailed above, $a_{1,2,3}$ have only variation $[-1,1]$ and $a_4$ has many orders of magnitude more in variation that $a_6$.
The conductor $N$ has a fairly tame distribution while the other BSD quantities vary rather wildly -- the is part of the difficulty of the conjecture, the relevant quantities behave quite unpredictably.

%%%%
\paragraph{The RHS of Conjecture 2: }
%%%
\begin{figure}[h!]
\centering
\includegraphics[width=\linewidth]{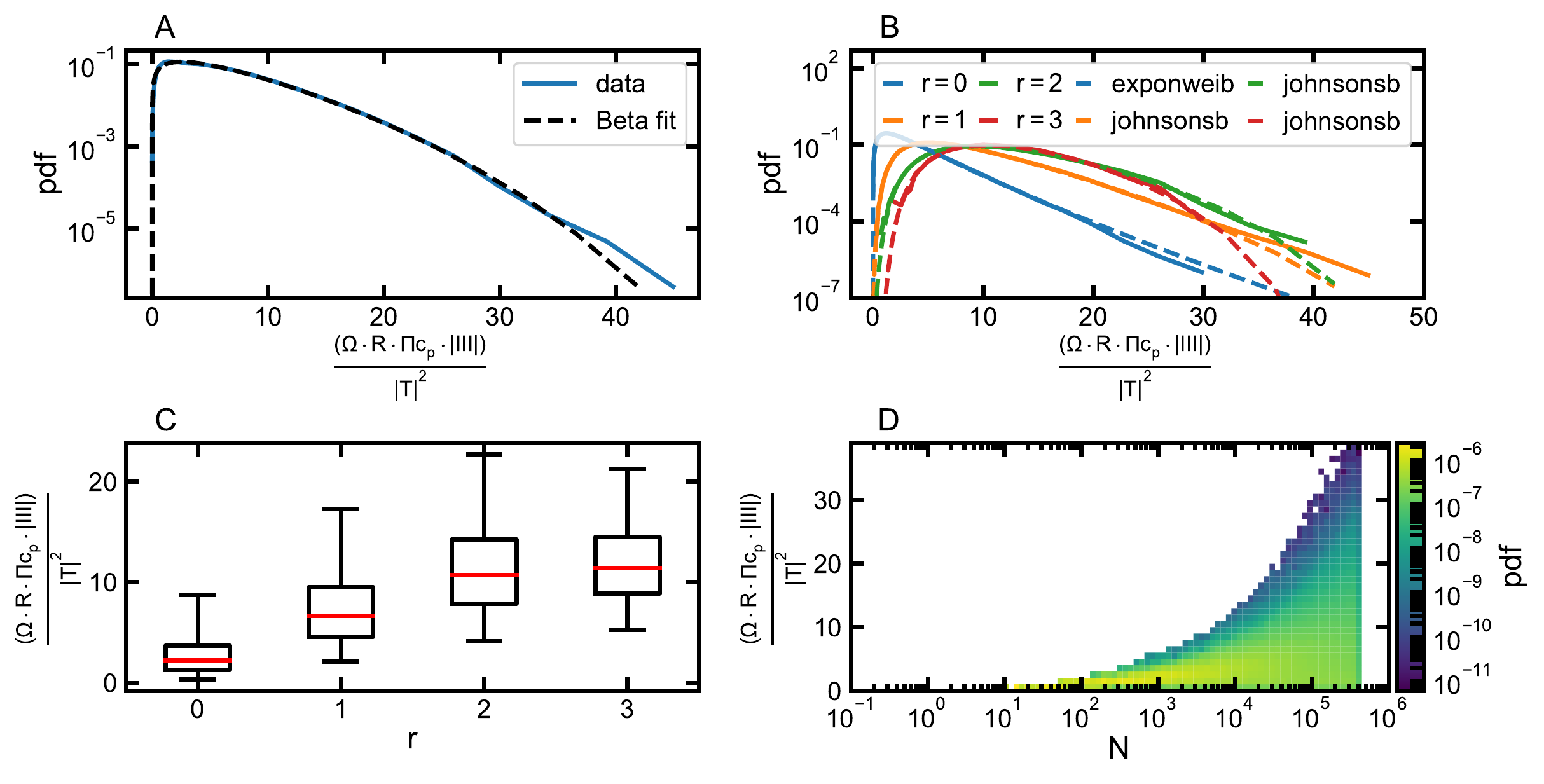}
\caption{\textbf{The RHS of Conjecture 2.} (A) Probability density function of the RHS (filled blue line) and the corresponding Beta distribution fit (dashed black line). (B) Probability density function of the RHS for different values of the rank $r$ (filled lines), and the corresponding best fits chosen with Akaike information criterion (dashed lines) (C) Boxplot showing the value of RHS for different values of $r$. The red line shows the median value, the boxes to the $50\%$ of the distribution, and the whiskers to the $95\%$ of the distribution.}
\label{RHS}
\end{figure}
%%%
We now put all the quantities together according to the RHS of the Strong BSD Conjecture, which we recall to be $RHS = \frac{(\Omega \cdot R \cdot \prod\limits_{p \mid N} c_p \cdot |\Sha|)}{T^2}$. 
We test which statistical distribution best describe this data, by comparing 85 continuous distribution under the Akaike information criterion. We find that the distribution best describing the data is the Beta distribution (see Figure \ref{RHS} A): 
\begin{equation}
f(x,a,b) = \frac{{\Gamma(a+b)x^{a-1}(1 - x)^{b-1}}}{\Gamma(a)\Gamma(b)}
\end{equation}
where $\Gamma$ is the standard gamma function, $a=1.55$, $b=14.28$,  and the original variable has been re-scaled such that $x = RHS/62.71$. 

The selected distribution changes for elliptic with a specific rank $r$.
For $r=0$, the selected distribution is the exponentiated Weibull, while for larger values of $r$, the selected distribution is the Johnson SB (see Figure \ref{RHS} B). 
We find that the median value of the $RHS$ increases both as a function of the rank $r$ (see Figure \ref{RHS} C) and $N$ (see Figure \ref{RHS} D).

%%%%%%%%%%%%%%%

%%%%%%%%%%%%%%%%%%%%%%%%
%%%%%%%%%%%%===================================
%%%%%%%%%%%%%%%%%%%%%%%%%%%%%
\section{Topological Data Analysis}
Let us gain some further intuition by visualizing the data.
As far as the data is concerned, to each elliptic curve, specified by the Weierstra\ss\  coefficients, one associates a list of quantities, the conductor, the rank, the real-period, etc.
The define a point cloud in Euclidean space of rather high dimension, each of point of which is defined by the list of these quantities which enter BSD.
In the above, we have extracted the last two coefficients of the Weierstra\ss\  form of the elliptic curves and studied them against the variations of the first three which, in normal form, can only be one of the 9 possible 3-tuples of $\pm1,$ and 0. The normal form thus conveniently allows us to at least ``see'' the Weierstra\ss\  coefficients because we have projected to two dimensions. However, the full list of the relevant quantities for the elliptic curve has quite a number of entries and cannot be visualized directly.

Luckily, there is precisely a recently developed method in data science which allows for the visualization of ``high dimensionality'', viz., persistent homolgy in topological data analysis \cite{czcg} (cf.~an excellent introductory survey of \cite{optgh}).
In brief, one creates a Vietoris-Rips simplex from the data points in Euclidean space, with a notion of neighbourhood $\epsilon$ (by Euclidean distance).
The Betti numbers $b_i$ of the simplex is then computed as one varies $\epsilon$, whether the values are non-zero for each $i$ gives a notion of whether non-trivial topology (such as holes) persists for different scales $\epsilon$.
The result is a so-called {\bf barcode} for the data.
In the ensuing, we will use G.~Henselman's nice implementation of the standard methods in topological data analysis, the package Eirene for Julia/Python \cite{eirene}.

%%%
\begin{figure}[h!]
\begin{tabular}{cc}
(a)
$b_0$ \includegraphics[width=2.5in]{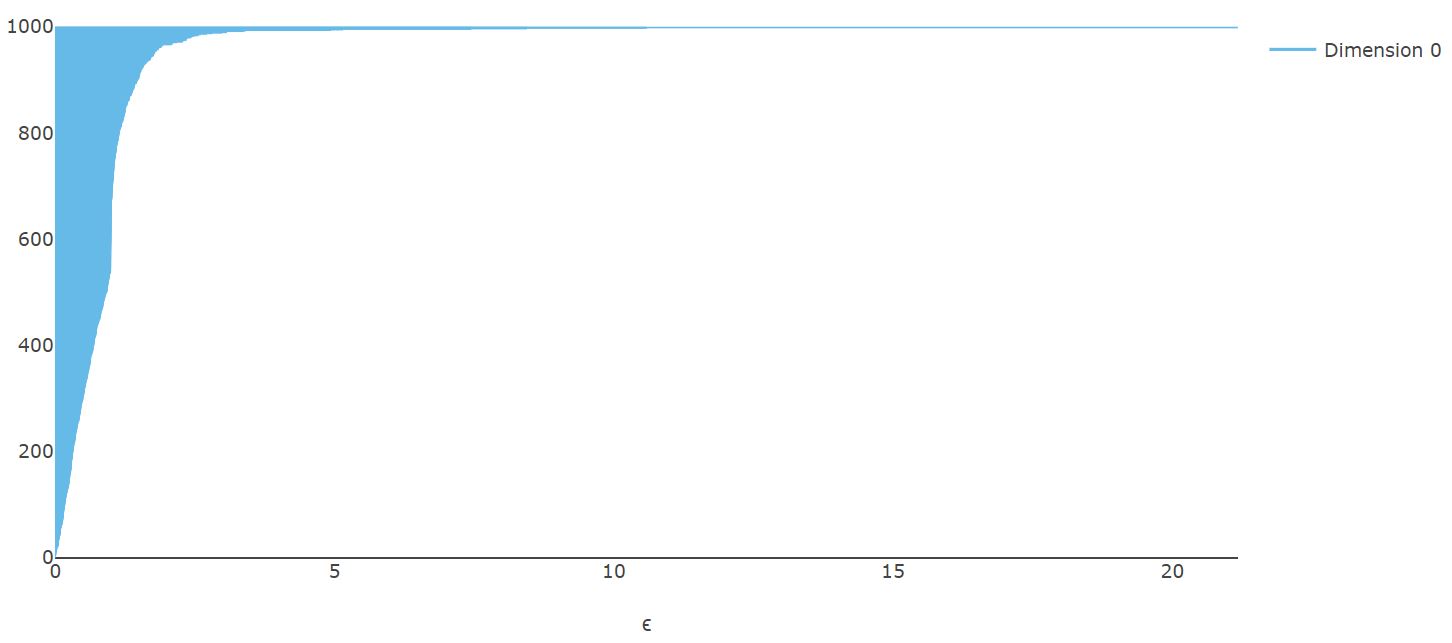} & $b_1$ \includegraphics[width=2.5in]{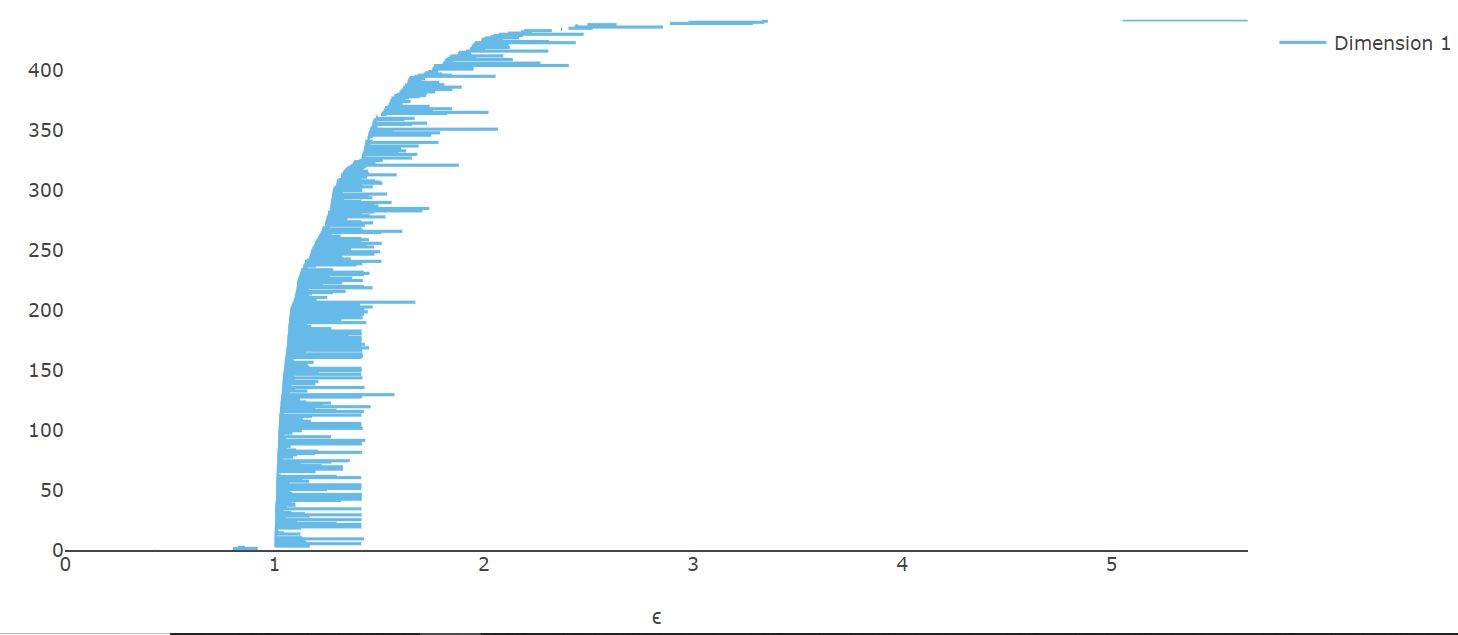} \\
(b)
$b_0$ \includegraphics[width=2.5in]{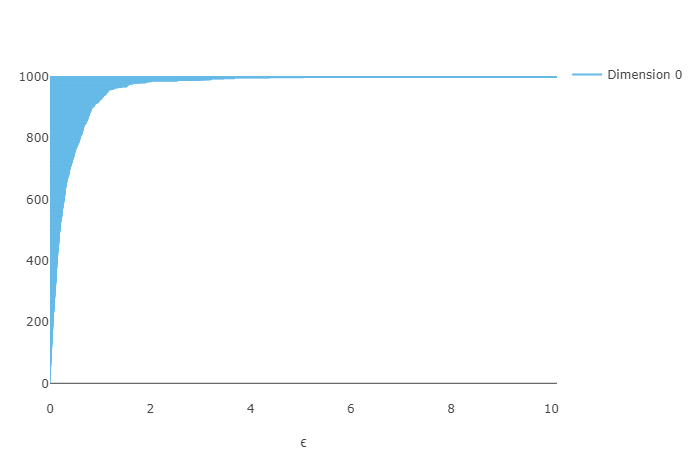} & $b_1$ \includegraphics[width=2.5in]{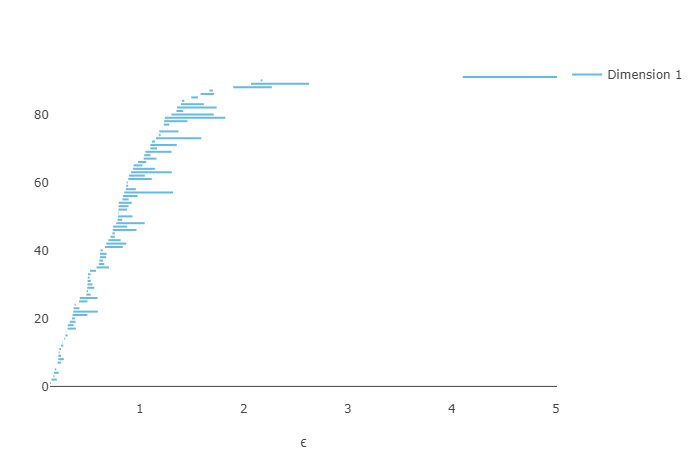} \\
\end{tabular}
\caption{{\sf The barcodes for (a) all the Weierstra\ss\  coefficients at Betti numbers 0 and 1; (b) on the (principle component) coefficients $\left( \slog(a_4), \ \slog(a_6) \right)$.}}
\label{f:coefbarcode}
\end{figure}
%%%

%%%%%%%%%%%%%%%%%%%%%
\subsection{The Weierstra\ss\  Coefficients}
We begin by studying the barcodes for the full set of five Weierstra\ss\  coefficients $a_i$ (as always, we take $\slog$ for $a_4$ and $a_6$ due to their size).
Of course, computing the homology of the Vietoris-Rips complex for over 2 million points is computationally impossible.
The standard method is to consider random samples.
Moreover, the first two Betti numbers $b_0$ and $b_1$ are usually sufficiently illustrative.
Thus, we will take 1000 random samples (we do not separate the ranks since we find there is no significant difference for the barcodes amongst different ranks) of the coefficients (with the usual $\slog$ for the last two).
The barcodes are shown in part (a) of  Figure \ref{f:coefbarcode}.
For reference, we also plot the barcodes for the pair $\left( \slog(a_4), \ \slog(a_6) \right)$ only since $a_{1,2,3}$ do not vary so much.

%%%%%%%%%%%%%%%%%%%%%
\subsection{A 6-dimensional Point-Cloud}
\begin{figure}[h!]
\begin{tabular}{ccc}
$r = 0$ & $r = 1$ & $r = 2$ \\
\includegraphics[width=2in]{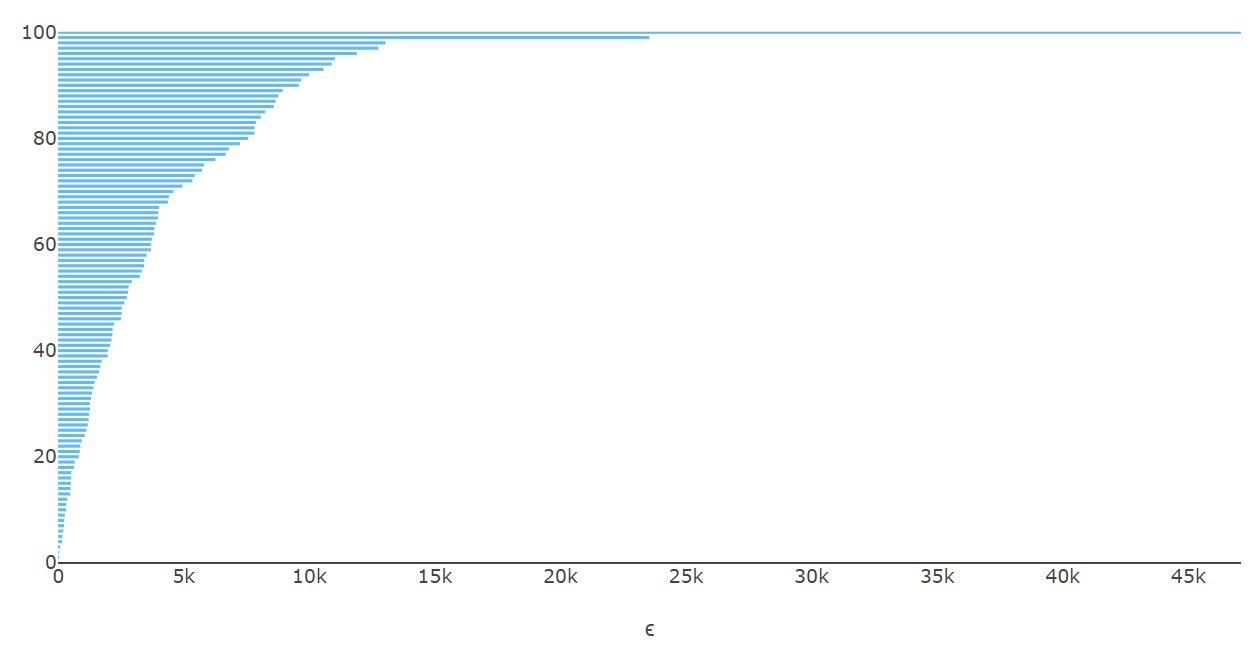} &
\includegraphics[width=2in]{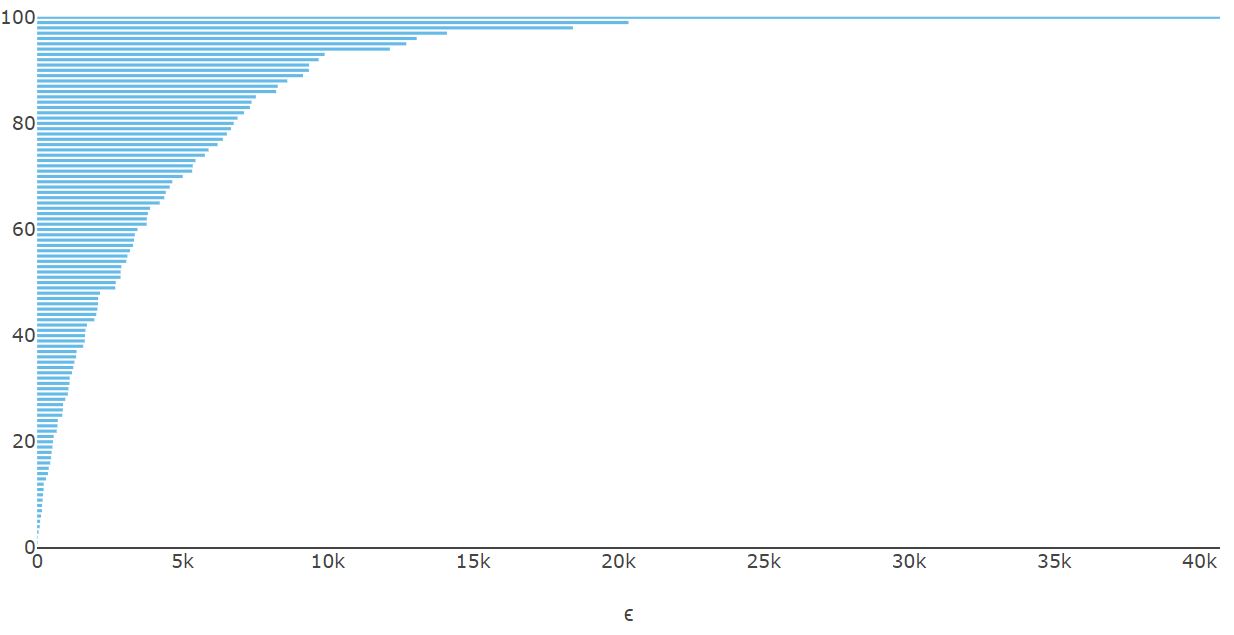} &
\includegraphics[width=2in]{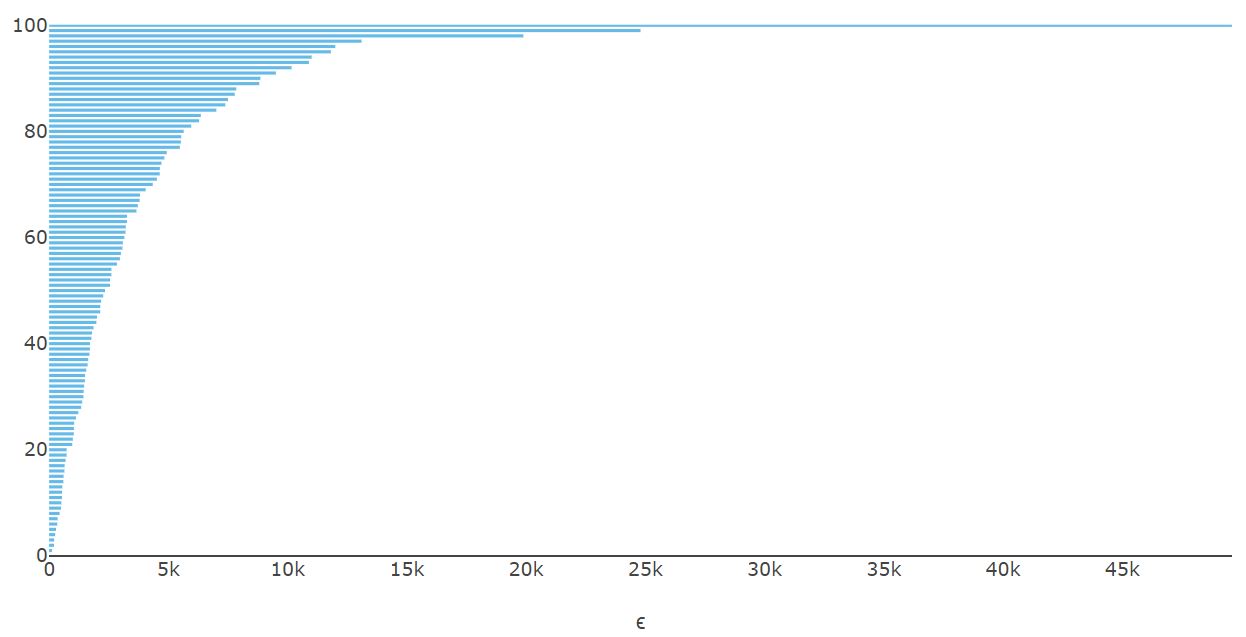}
\end{tabular}
\caption{{\sf The barcodes for 6-dimensional point-cloud $(N,\ |T|,\ \prod\limits_{p \mid N} c_p,\ \Omega,\ R, \Sha)$, with 100 random samples, for $r=0,1,2$.}}
\label{f:6barcode}
\end{figure}

Let us now try to visualize the relevant BSD quantities $(N,\ r,\ |T|,\ \prod\limits_{p \mid N} c_p,\ \Omega,\ R, \Sha)$ together.
Organized by the ranks $r = 0,1,2$ which dominate the data by far, the 6-tuple
\begin{equation}
(N,\ |T|,\ \prod\limits_{p \mid N} c_p,\ \Omega,\ R, \Sha), \qquad r =0,1,2 
\end{equation}
naturally form three point-clouds in $\IR^6$.
Due to the high dimensionality, we sample 100 random points for each of the $r$ values and compute the full barcodes $b_{0,\ldots,6}$.
It turns out that the main visible features are in dimension 0.
We present these in Figure \ref{f:6barcode} and observe that indeed there is some variation in the barcode amongst the different ranks.

%%%%%%%%%%%%%%%%%%%%%
\subsection{Conductor Divisibility}
\begin{figure}[h!]
\begin{tabular}{cc}
$N$ even &  $N$ odd \\
\includegraphics[width=3in]{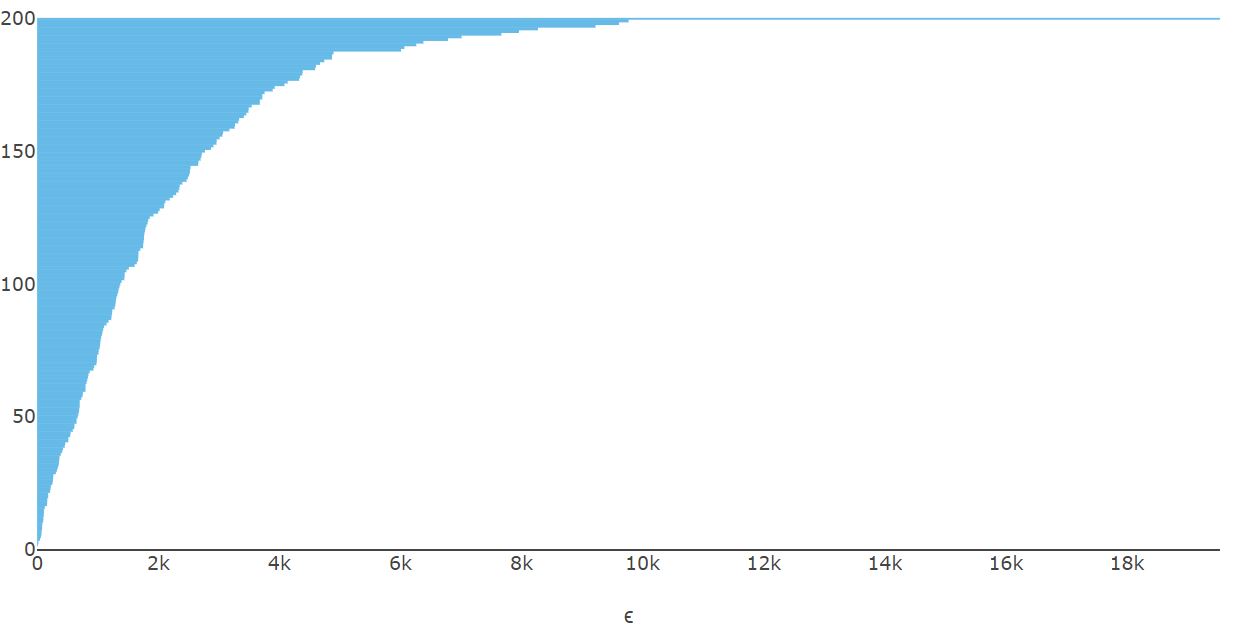} &
\includegraphics[width=3in]{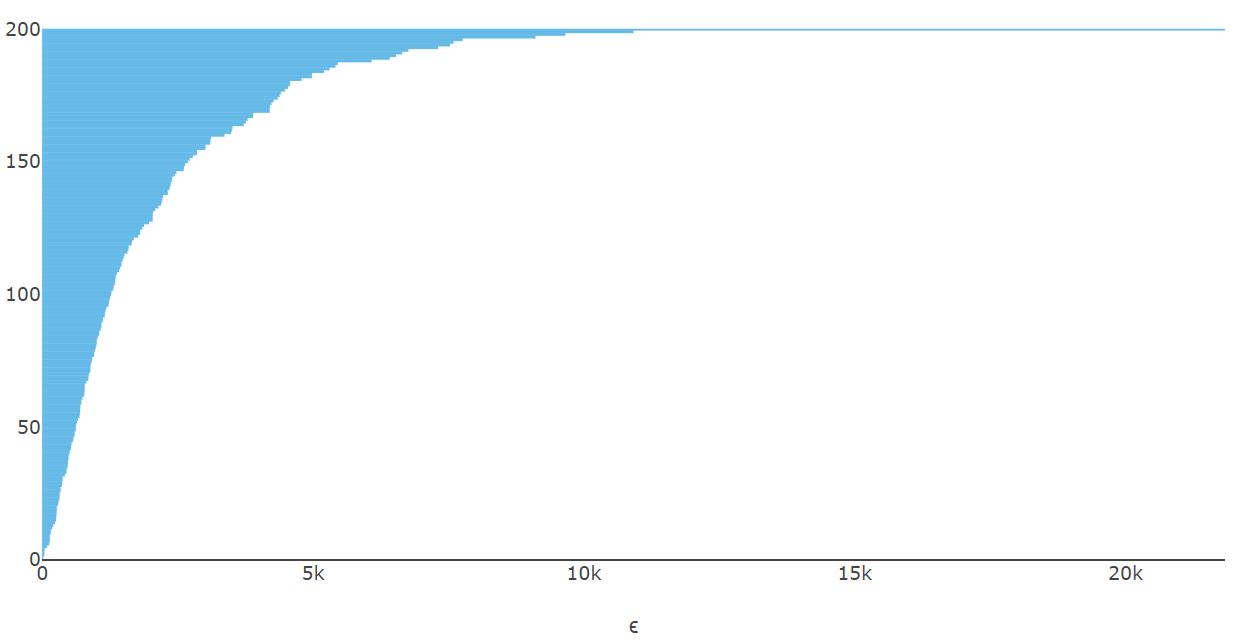}
\end{tabular}
\caption{{\sf The barcodes for 3-dimensional point-cloud $(N,\ r,\ \frac{|\Sha| \cdot \Omega \cdot R \cdot \prod\limits_{p \mid N} c_p}{|T|^2})$, with 200 random samples for even/odd $N$.}}
\label{f:Nmod2}
\end{figure}

The factors of the conductor $N$ is of defining importance in the L-function, which appear to the LHS of BSD, meanwhile, the RHS is governed, in the strong case, by the combination 
$F := \frac{|\Sha| \cdot \Omega \cdot R \cdot \prod\limits_{p \mid N} c_p}{|T|^2}$.
It is therefore expedient to consider the point cloud of the 3-tuple $(N, r, F)$ organized by divisibility properties of $N$.
For instance, one could contrast the barcodes for the triple for $N$ even versus $N$ odd.
Again, the features are prominent for dimension 0 and the barcodes are shown in Figure \ref{f:Nmod2}.

\begin{figure}[h!]
\begin{tabular}{ccc}
$N\equiv 0$ &  $N \equiv 1$ & $N\equiv 2$ \\
\includegraphics[width=2in]{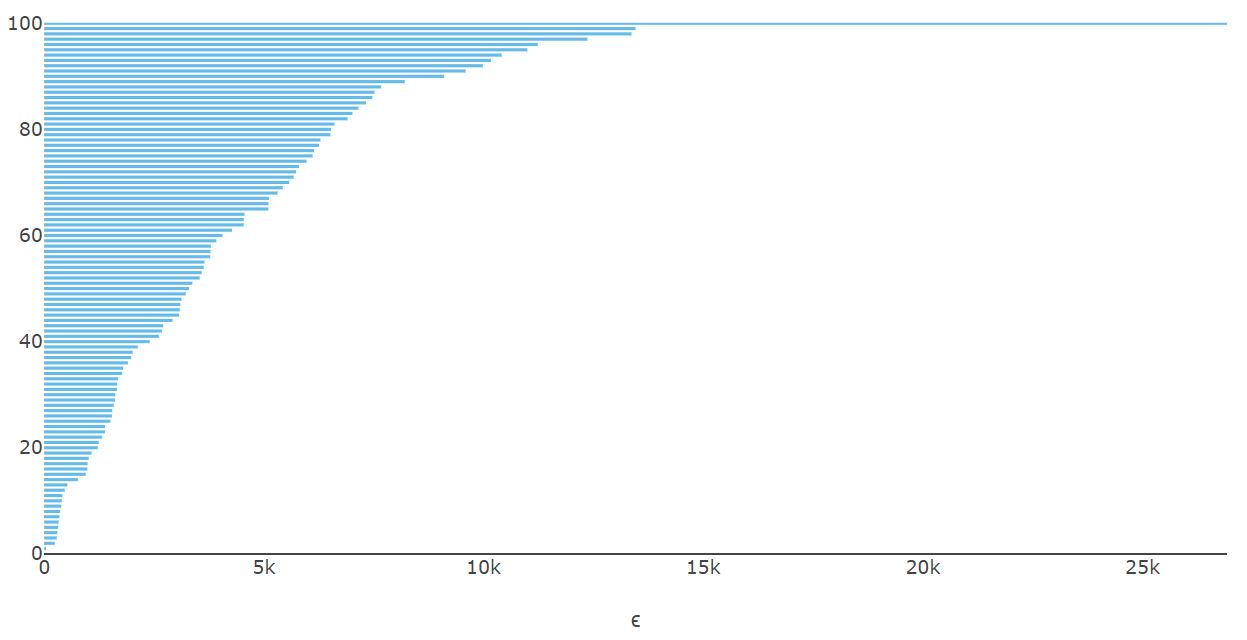} &
\includegraphics[width=2in]{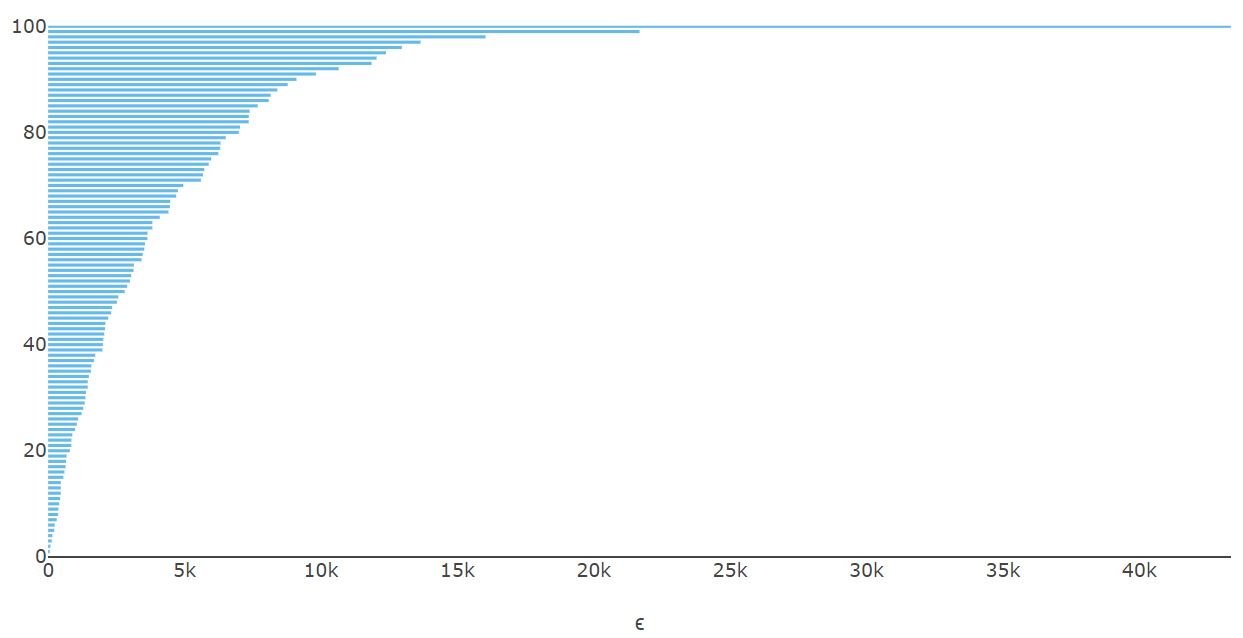} &
\includegraphics[width=2in]{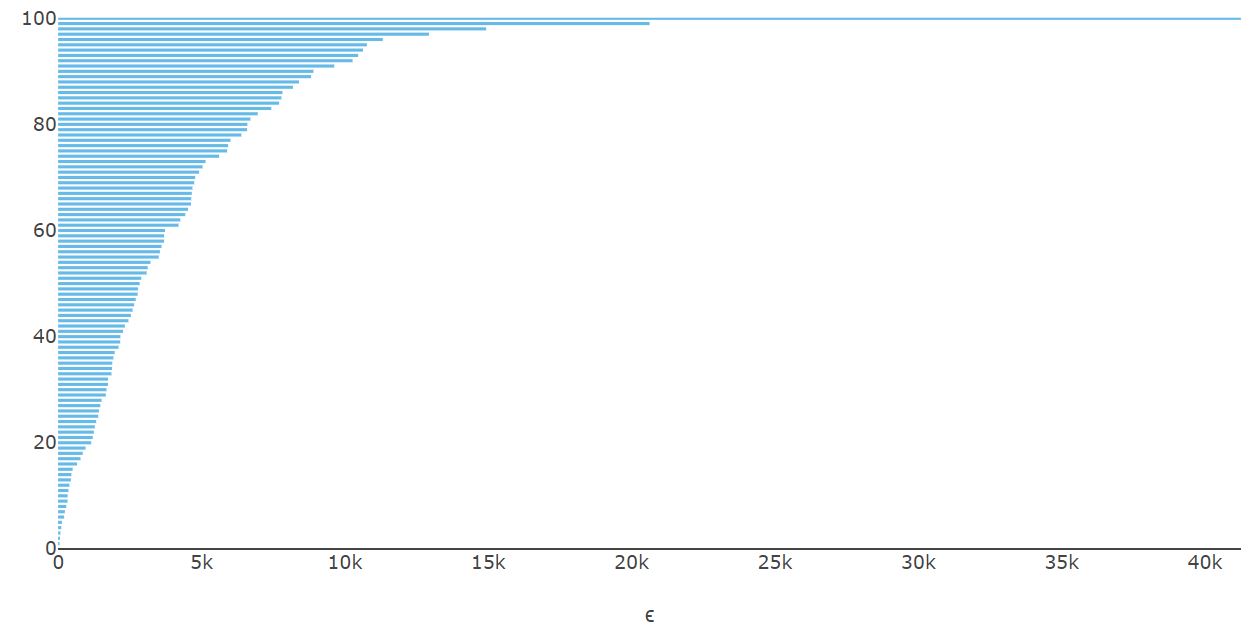}
\end{tabular}
\caption{{\sf The barcodes at dimension 0 for 3-dimensional point-cloud $(N,\ r,\ \frac{|\Sha| \cdot \Omega \cdot R \cdot \prod\limits_{p \mid N} c_p}{|T|^2})$, with 100 random samples, for $N$ distinguished modulo 3.}}
\label{f:Nmod3}
\end{figure}

Simiarly, we could group by $N$ modulo 3, as shown in Figure \ref{f:Nmod3}.

%%%%%%%%%%%%%%%
\section{Machine Learning}
In \cite{He:2017aed,He:2018jtw}, a paradigm was proposed to using machine-learning, and in particular deep neural networks to help computations in various problems in algebraic geometry.
Exemplified by computing cohomology of vector bundles, it was shown that to very high precision, AI can guess the correct answer without using the standard method of Gr\"obner basis construction and chasing long exact sequences, both of which are computationally intensive.
Likewise, \cite{He:2019nzx} showed that machine-learning can identify algebraic structures such as distinguishing simple from non-simple finite groups.
At over 99\% precision, the requisite answers can be estimated without recourse to standard computations which are many orders of magnitude slower.

It is therefore natural to wonder whether the elliptic curve data can be ``machine-learned''.
Of course, we need to be careful.  While computational algebraic geometry over $\IC$ hinged on finding kernels and cokernels of integer matrices, a task in which AI excels. Problems in number theory are much less controlled. Indeed, trying to predict prime numbers \cite{He:2017aed} seems like a hopeless task, as mentioned in the introduction.
Nevertheless, let us see how far we can go with our present dataset for BSD.

%%%%%%%%%%%%%%%%%%%%%
\subsection{Predicting from the Weierstra\ss\ Coefficients}\label{ml1}

%%%% 
\begin{figure}[h!]
\centering
\includegraphics[width=.7\linewidth]{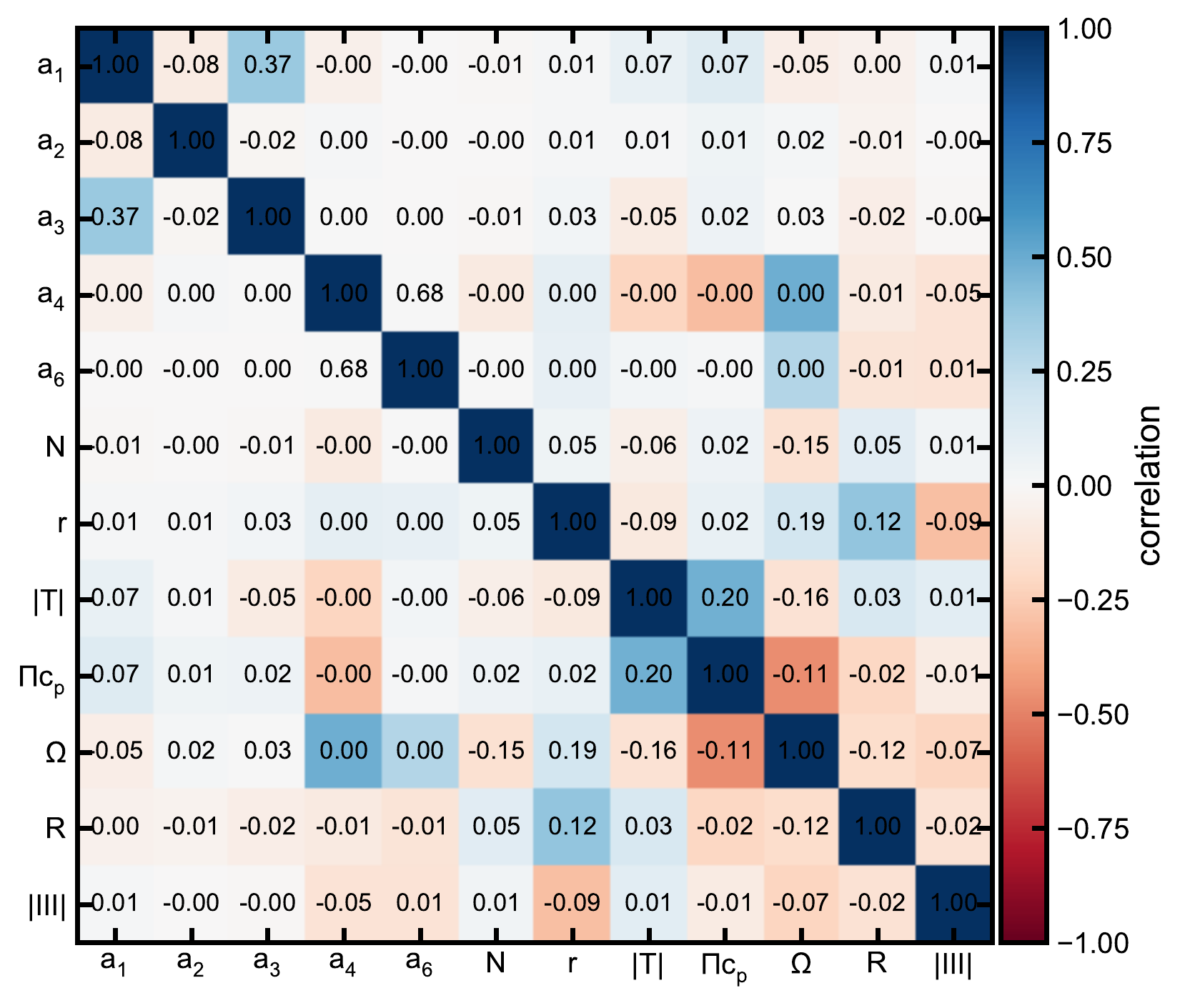}
\caption{\textbf{Correlation matrix for the quantities characterizing ellipses.} Colours are assigned based on the value of the Pearson correlation coefficient between all pairs of quantities. The strength of the correlation is also reported for each pair. }
\label{correlation_matrix}
\end{figure}
%%%%%

We begin with quantifying the performance of machine learning models in predicting, one by one, the quantities $N$, $\prod\limits_{p \mid N} c_p$, $R$, $|\Sha|$, $\Omega$, $r$, $|T|$, together with the RHS of the Conjecture 2, given the Weierstra\ss\ coefficients $a_1$, $a_2$, $a_3$, $a_4$, $a_6$ alone. 
Straight-away, this is expected to be a difficult, if not impossible task (as impossible as, perhaps, the prediction of prime numbers). This is confirmed by the correlation between the Weierstra\ss coefficients and the BSD quantities: from the correlation matrix, we see that the relationship is indeed weak (cf.~Fig.~\ref{correlation_matrix}), implying this is not a straightforward prediction task. Here, the correlation matrix contains the values of the Spearman correlation \cite{pearson}, that is the Pearson correlation \cite{pearson} between the rank variables associated to each pairs of the original variables.

%%%
\begin{table}[!ht]
    \centering
    \begin{tabularx}{\columnwidth}{@{}Xrrr@{}}
        \midrule
        \textbf{Quantity} & \textbf{NMAE (XGBoost)} & \textbf{NMAE (Dummy)} & \textbf{NMAE (Linear)} \\ \toprule
$N$ & $24.842\pm 0.032 $& $25.175\pm 0.026$& $25.158\pm 0.027$ \\ 
$\prod\limits_{p \mid N} c_p$ & $0.028\pm 0.006 $& $0.077\pm 0.016$& $0.058\pm 0.012$ \\ 
$R$ & $0.075\pm 0.015 $& $0.112\pm 0.023$& $0.108\pm 0.022$ \\ 
$|\Sha|$ & $0.023\pm 0.015 $& $0.044\pm 0.028$& $0.043\pm 0.027$ \\ 
$\Omega$ & $3.120\pm 0.099 $& $6.057\pm 0.189$& $6.016\pm 0.189$ \\ 
RHS Conj. 2 & $7.070\pm 0.238 $& $7.548\pm 0.255$& $7.533\pm 0.250$ \\
         \bottomrule
    \end{tabularx}
   \begin{tabularx}{\columnwidth}{@{}Xrrr@{}}
        \midrule
         & \textbf{RMSE (XGBoost)} & \textbf{RMSE (Dummy)} & \textbf{RMSE (Linear)} \\ \toprule
$N$ & $114687.179\pm 63.171 $& $115784.768\pm 78.329$& $115774.283\pm 78.302$ \\ 
$\prod\limits_{p \mid N} c_p$ & $273.912\pm 18.665 $& $286.522\pm 19.679$& $285.731\pm 19.711$ \\ 
$R$ & $13.579\pm 0.886 $& $14.201\pm 0.552$& $14.197\pm 0.555$ \\ 
$|\Sha|$ & $6.797\pm 1.550 $& $6.369\pm 1.794$& $6.524\pm 1.688$ \\ 
$\Omega$ & $0.449\pm 0.001 $& $0.584\pm 0.001$& $0.583\pm 0.001$ \\ 
RHS Conj. 2 & $4.300\pm 0.002 $& $4.554\pm 0.004$& $4.526\pm 0.003$ \\
         \bottomrule
    \end{tabularx}
    \caption{\textbf{Performance of the regression models.} The Normalized Median Absolute Error ($NMAE$) and the Root Mean Squared Error ($RMSE$), for XGboost (left column), the dummy regressor (central column) and a linear regression (right column). The reported values are averages across $5$-fold cross-validations, with the corresponding standard deviations.}
\label{results_regression_1}
\end{table}
%%%%

Our analysis relies on gradient boosted trees \cite{friedman2001greedy}, using the implementation of XGboost \cite{chen2016xgboost}, an open-source scalable machine learning system for tree boosting used in a number of winning Kaggle solutions (17/29 in 2015). 
In Appendix \ref{svm}, we present the similar results using a support vector machine, another highly popular machine-learning model, and see that the XGboost indeed performs better.
Furthermore, based on the learning curves of the XGboost models (discussed in Appendix \ref{learning_curve}), we have chosen a {\bf $5-$fold cross-validation}, such that the training set includes $80\%$ of the values, and the validation set the remaining $20\%$.

\subsubsection{Numerical Quantities}
First, we train regression models to predict the values of $N$, $\prod\limits_{p \mid N} c_p$, $R$, $|\Sha|$ and $\Omega$. We evaluate the performance of the regression, by computing the normalized median absolute error: 
\begin{equation}
NMAE = \frac{median(|Y_{i} - \hat{Y}_{i}|)}{max(Y_{i}) - min(Y_{i})} \ ,
\end{equation}
where $Y_i$ are the observed values and $\hat{Y}_{i}$ are the predicted values, and the rooted mean squared error: 
\begin{equation}
RMSE = \sqrt{\frac{\sum(Y_{i} - \hat{Y}_{i})^2}{n}} \ ,
\end{equation}
where $n$ is the size of the test set. 
We desire that both NMAE and RMSE to be close to 0 for a good prediction.

We compare the result of the XGBoost regression with two baselines: (1) a linear regression model and (2) a dummy regressor, that always predicts the mean of the training set. 
We find that, in all cases, the machine learning algorithms perform significantly better than the baseline models (see Table \ref{results_regression_1}) with respect to the NMAE and RMSE. However, XGboost performs only marginally better than the baselines in predicting the value of $N$. We report also the so-called {\im importance} of the features for the XGBoost regressor in Fig.~\ref{regr_1}.  
Here, importance indicates how useful each feature is in the construction of the boosted decision trees, and is calculated as the average importance across trees. For a single tree, the importance of a feature is computed as the relative increase in performance resulting from the tree splits based on that given feature \cite{chen2016xgboost}.

%%%%%%%%
\begin{figure}[h!]
\centering
\includegraphics[width=\linewidth]{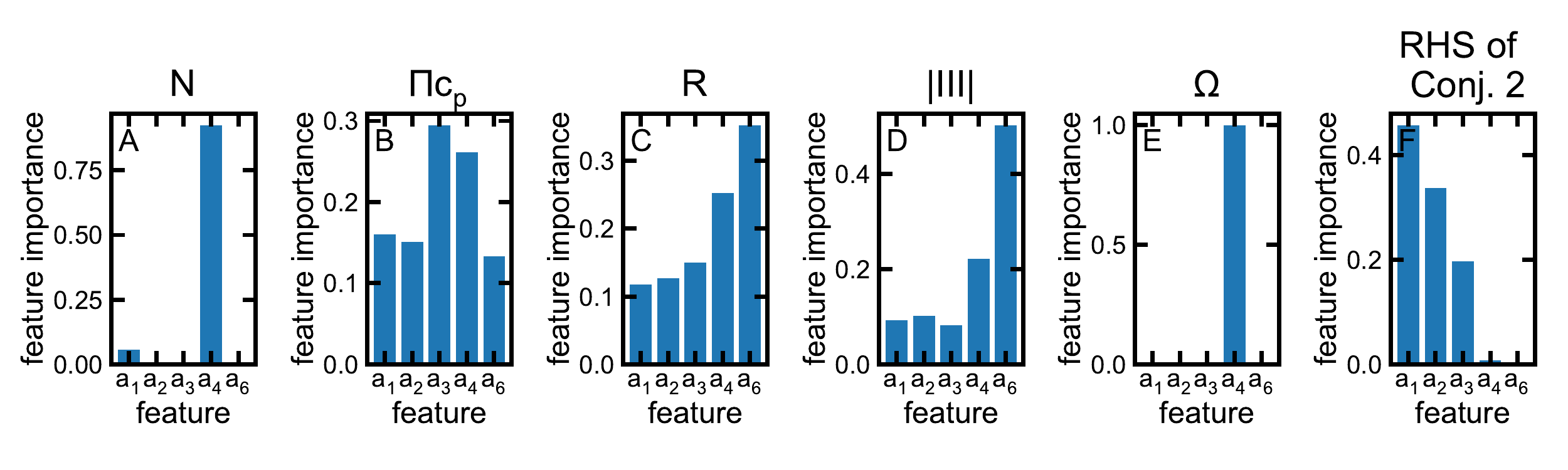}
\caption{\textbf{Feature importance.} Feature importance of the XGBoost regression models for predicting $N$ (A), $\prod\limits_{p \mid N} c_p$ (B), $R$ (C), $|\Sha|$ (D) and $\Omega$ (E).}
\label{regr_1}
\end{figure}
%%%%%%%%

We see that overall our measures, NMAE and RMSE, are not too close to 0, except perhaps $|\Sha|$, for which even a simple linear regression does fairly well. 
In table \ref{lin1}, we report the values of the coefficients of the linear regression fit for $|\Sha|$.
\begin{table}[h!!!]
\begin{center}
\begin{tabular}{lrrrrrr}
\toprule
{} &    coef &  std err &        t &  P>|t| &  [0.025 &  0.975] \\
\midrule
const &  1.5946 &    0.004 &  380.761 &  0.000 &   1.586 &   1.603 \\
$a_1$    &  0.0658 &    0.005 &   14.524 &  0.000 &   0.057 &   0.075 \\
$a_2$    & -0.0065 &    0.004 &   -1.543 &  0.123 &  -0.015 &   0.002 \\
$a_3$    & -0.0518 &    0.005 &  -11.473 &  0.000 &  -0.061 &  -0.043 \\
$a_4$   & -0.6320 &    0.006 & -110.282 &  0.000 &  -0.643 &  -0.621 \\
$a_6$    &  0.4877 &    0.006 &   85.112 &  0.000 &   0.477 &   0.499 \\
\bottomrule
\end{tabular}
\end{center}
\caption{\textbf{Prediction of $|\Sha|$.} Coefficients of the linear regression for each of the features (inputs), with the associated standard deviation, the value of the $t$ statistics and corresponding $p$-value. Here, we can reject the null hypothesis that the coefficients are equal to $0$ at significance level $\alpha=0.01$, for all coefficients, except the one associated to $a_2$ ($p>0.01$).}
\label{lin1}
\end{table}
%%%
What Table \ref{lin1} means is that $|\Sha| \simeq 1.5946 + 0.0658 a_1 -0.0065 a_2  -0.0518 a_3 - 0.6320  a_4 + 0.4877a_6$.

\begin{table}[h!!!]
\begin{center}
\begin{tabular}{rrrrr}
\toprule
Predicted variable & R-squared & Adj. R-squared & F-statistic & Prob (F-statistic) \\
$N$ & 0 &               0 &        95.67 &           3.86e-101 \\
$\prod\limits_{p \mid N} c_p$ &    0.006 &           0.006 &         2777 &                   0 \\
$R$ &      0.001 &           0.001 &        387.2 &                   0  \\
$|\Sha|$ &     0.005 &           0.005 &         2522 &                   0  \\
$\Omega$ & 0.005 &           0.005 &         2377 &                   0 \\
RHS Conj. 2 &     0.012 &           0.012 &         6143 &                   0 \\
\bottomrule
\end{tabular}
\end{center}
\caption{\textbf{Statistics of the linear regression models.} We report the R squared, the adjusted R squared, the F statistics, and the p-value associated to the F statistics for the different linear models. When the p-value of the F-statistics is close to $0$, we can reject the null hypothesis  that the intercept-only model provides a better fit than the linear model. }
\label{lin2}
\end{table}
%%%

Likewise, in table \ref{lin2}, we report the the statistics associated to the linear regression models for the prediction of the various quantities, where only the Weierstrass coefficients are used as features (inputs).
The low R-squared values indicated that the regression is not so good in terms of the $a_i$ coefficients alone.

%%%%%%%%%%%%%%%%%%%%%%%%
\subsubsection{Categorical Quantities}
Next, we train classifiers to predict the values of $r$ and $|T|$ because these easily fall into discrete categories:
the rank $r = 0,1,2,3$ and the torsion group size $|T|$ can only be one of 16 integer values due to Mazur's theorem.
Again, we use a $5-fold$ cross validation, and we evaluate the performance of the classifier, by computing the $F1$ score:
\begin{equation}
F1 = 2 \cdot \frac{\mbox{precision} \cdot \mbox{recall}}{\mbox{precision}+\mbox{recall}} \ ; \quad
\mbox{precision} := \frac{TP}{TP + FP} \ , \ 
\mbox{recall} := \frac{TP}{TP + FN}
\end{equation}
where we have, in the predicted versus actual, the true positives (TP), false positives (FP), and false negatives (FN).
Since we have several possible values for the rank $r$, we compute both $F1_{micro}$, by counting the total TP, FN and FP, as well as $F1_{macro}$, the average value of $F1$ computed across ranks. 

In addition, we also compute the Matthew correlation coefficient $MCC$ \cite{matthew}, to describe the confusion matrix:
\begin{equation}
MCC := \frac{ \mathit{TP} \times \mathit{TN} - \mathit{FP} \times \mathit{FN} } {\sqrt{ (\mathit{TP} + \mathit{FP}) ( \mathit{TP} + \mathit{FN} ) ( \mathit{TN} + \mathit{FP} ) ( \mathit{TN} + \mathit{FN} ) } } \ .
\end{equation}
Both F1-scaore and MCC are desired to be close to 1 for a good prediction.

%%%%%%
\begin{table}[ht]
    \centering
    \begin{tabularx}{\columnwidth}{@{}Xrrr@{}}
        \midrule
        \textbf{Quantity} & \textbf{ $F1_{micro}$ (XGBoost)} & \textbf{ $F1_{micro}$ (Dummy)} & \textbf{ $F1_{micro}$ (Logistic)} \\ \toprule
$r$ & $0.502\pm 0.001 $& $0.502\pm 0.001$& $0.502\pm 0.001$ \\ 
$|T|$ & $0.582\pm 0.001 $& $0.543\pm 0.001$& $0.518\pm 0.001$ \\ 
         \bottomrule
    \end{tabularx}
    \begin{tabularx}{\columnwidth}{@{}Xrrr@{}}
        \midrule & \textbf{ $F1_{macro}$ (XGBoost)} & \textbf{ $F1_{macro}$ (Dummy)} & \textbf{ $F1_{macro}$ (Logistic)} \\ \toprule
$r$ & $0.179\pm 0.001 $& $0.167\pm 0.001$& $0.167\pm 0.001$ \\ 
$|T|$ & $0.097\pm 0.001 $& $0.059\pm 0.001$& $0.080\pm 0.001$ \\ 
         \bottomrule
    \end{tabularx}
    \begin{tabularx}{\columnwidth}{@{}Xrrr@{}}
        \midrule & \textbf{ $MCC$ (XGBoost)} & \textbf{ $MCC$ (Dummy)} & \textbf{ $MCC$ (Logistic)} \\ \toprule
$r$ & $0.0172\pm 0.0006 $& $0.0000\pm 0.0000$& $-0.0002\pm 0.0010$ \\ 
$|T|$ & $0.1871\pm 0.0010 $& $0.0000\pm 0.0000$& $0.0299\pm 0.0012$ \\
         \bottomrule
    \end{tabularx}
    \caption{\textbf{Performance of the classification models.} The scores $F1_{micro}$,  $F1_{macro}$, and the Matthew correlation coefficient $MCC$, for XGboost (left column), the dummy regressor (central column) and a logistic regression (right column). The reported values are averages across $5$-fold cross-validations, with the corresponding standard deviations. }
\label{results_classification}
\end{table}
%%%%%%

For checks, we compare the XGBoost classifier with (1) a baseline classifier, that predicts always the predominant class in the training set, as well as with (2) a logistic regression. We find that XGboost performs better than the baseline models for predicting $|T|$, but the performance of the prediction of $r$  is comparable to the baselines (see Table \ref{results_classification}).

%%%%%%%%%%%%
\begin{figure}[h!!!]
\centering
\includegraphics[width=\linewidth]{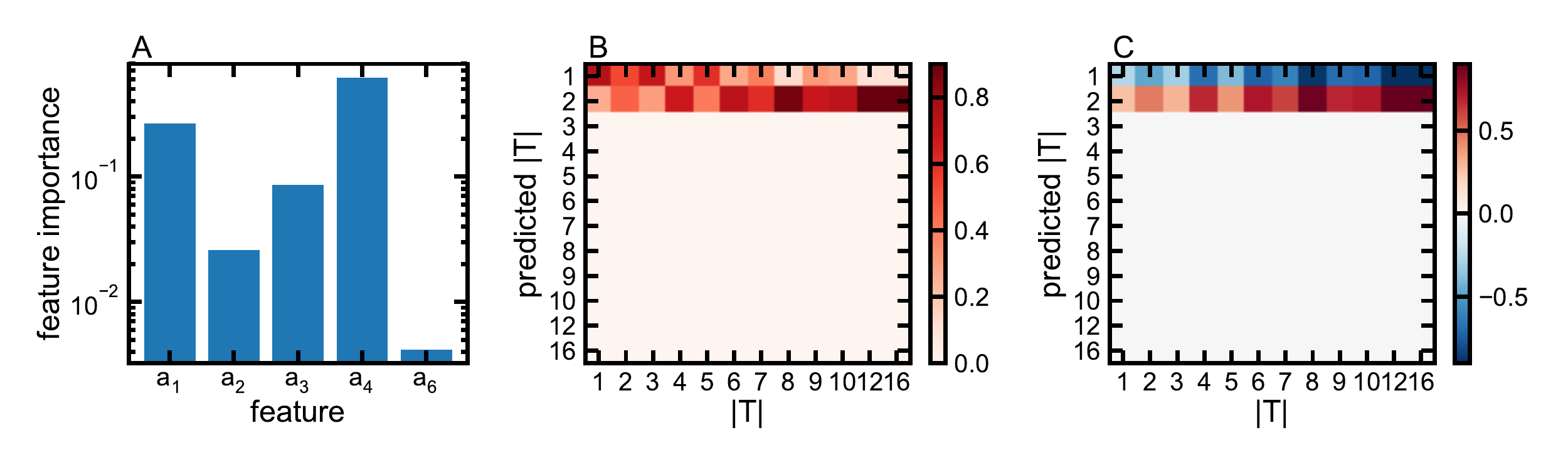}
\caption{\textbf{Prediction of $|T|$.} (A) Importance of the different features (inputs) to predict $|T|$. (B) Confusion matrix (normalized by column) showing the fraction of entries $|T|$ with given $predicted$ $|T|$. (C) Difference between the confusion matrix obtained for the XGBoost and the dummy classifier. Results are averaged over a $5$-fold cross validation.}
\label{class_1}
\end{figure}

\begin{figure}[h!!!]
\centering
\includegraphics[width=\linewidth]{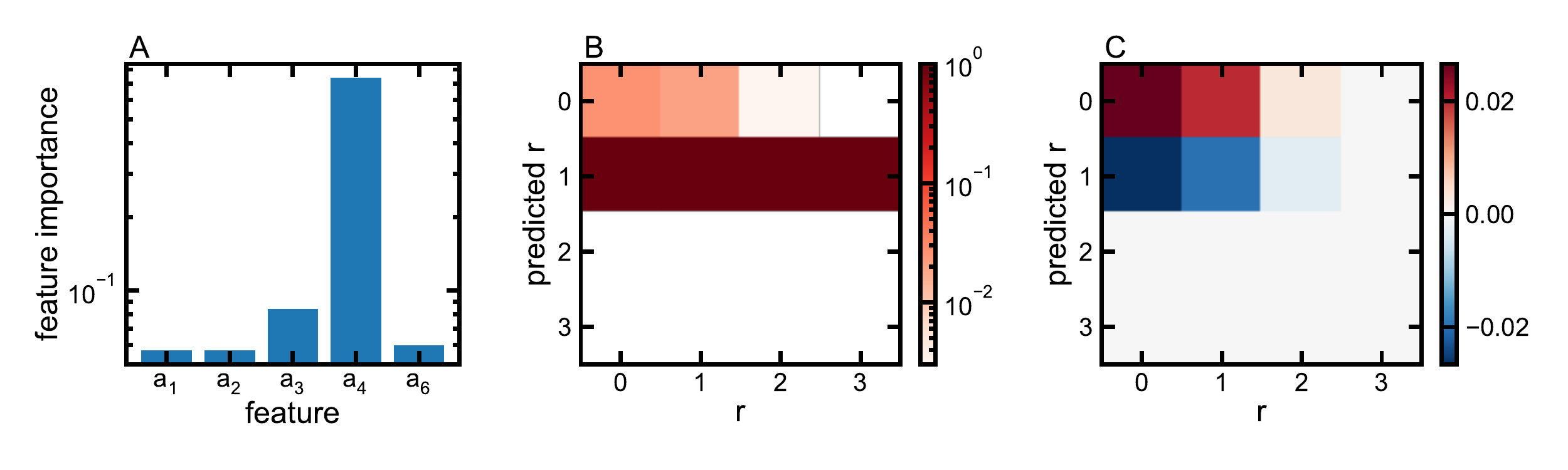}
\caption{\textbf{Prediction of $r$.} (A) Importance of the different features to predict $r$. (B) Confusion matrix (normalized by column) showing the fraction of entries with rank $r$ with given $predicted$ $r$. (C) Difference between the confusion matrix obtained for the XGBoost and the dummy classifier. Results are averaged over a $5$-fold cross validation.}
\label{class_2}
\end{figure}
%%%%%%%%%%%

Analyzing the confusion matrices (see Figures \ref{class_1} and \ref{class_2}), it appears clear that it is very hard to predict $|T|$ and $r$ from the Weierstra\ss\ coefficients alone. For both the prediction of $r$ and $|T|$, the most important predictor is $a_4$ (Figures \ref{class_1} and \ref{class_2}). This is the feature that contributed the most to increase the performance of the boosted tree \cite{chen2016xgboost}.

%%%%%%%%%%%%%%%%%%%%%%%%%%%=====================
%%%%%%%%%%%%%%%%%%%%%%%%%
\subsection{Mixed Predictions}
While the results in the previous subsection may seem disappointing in general, they do present a good sanity check: to obtain all the BSD quantities from the elliptic curve data in some straight-forward way would be an almost unimaginable feat in number theory.
Nevertheless, in this section, let us build machine learning models to predict the values of $N$, $\prod\limits_{p \mid N} c_p$, $R$, $|\Sha|$, $\Omega$, $r$ and $|T|$ among themselves, i.e., we consider as features (inputs) all the quantities characterizing the elliptic curves (except the predicted quantity), rather than the Weierstra\ss\ coefficients alone. 

%%%
\begin{table}[h!!!]
    \centering
    \begin{tabularx}{\columnwidth}{@{}Xrrr@{}}
        \midrule
        \textbf{Quantity} & \textbf{NMAE (XGBoost)} & \textbf{NMAE (Dummy)} & \textbf{NMAE (Linear)} \\ \toprule
$N$ & $23.426\pm 0.031 $& $25.175\pm 0.026$& $24.408\pm 0.039$ \\ 
$\prod\limits_{p \mid N} c_p$ & $0.012\pm 0.003 $& $0.077\pm 0.016$& $0.065\pm 0.014$ \\ 
$R$ & $0.014\pm 0.003 $& $0.112\pm 0.023$& $0.089\pm 0.018$ \\ 
$|\Sha|$ & $0.006\pm 0.004 $& $0.044\pm 0.028$& $0.048\pm 0.031$ \\ 
$\Omega$ & $2.343\pm 0.103 $& $6.057\pm 0.189$& $5.324\pm 0.174$ \\ 

         \bottomrule
    \end{tabularx}
    \caption{\textbf{Performance of the regression models considering as features all the quantities characterizing an ellipses.} The normalized median absolute error $NMAE$, for XGboost (left column), the dummy regressor (central column) and a linear regression (right column). The reported values are averages across $5$-fold cross-validations, with the corresponding standard deviations. Results are considerably improved compared to \cref{results_regression_1}.}
\label{results_regression_3}
\end{table}

\begin{figure}[h!]
\centering
\includegraphics[width=\linewidth]{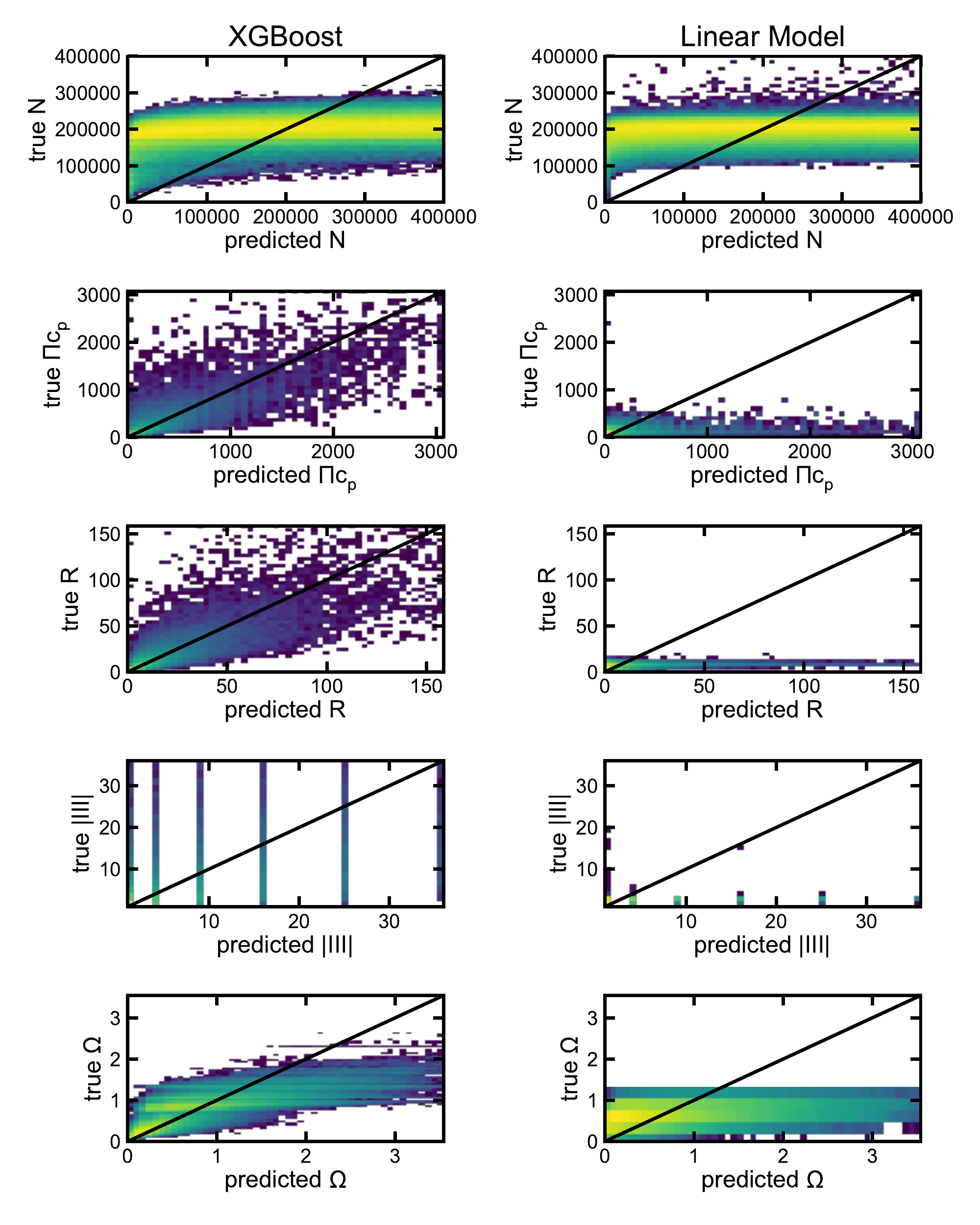}
\caption{\textbf{True vs Predicted values} Results are shown for all the quantities, using XGBoost (left column) and the Linear model (right column).}
\label{Figure_True_pred}
\end{figure}

%%%

\begin{figure}[h!!!]
\centering
\includegraphics[width=\linewidth]{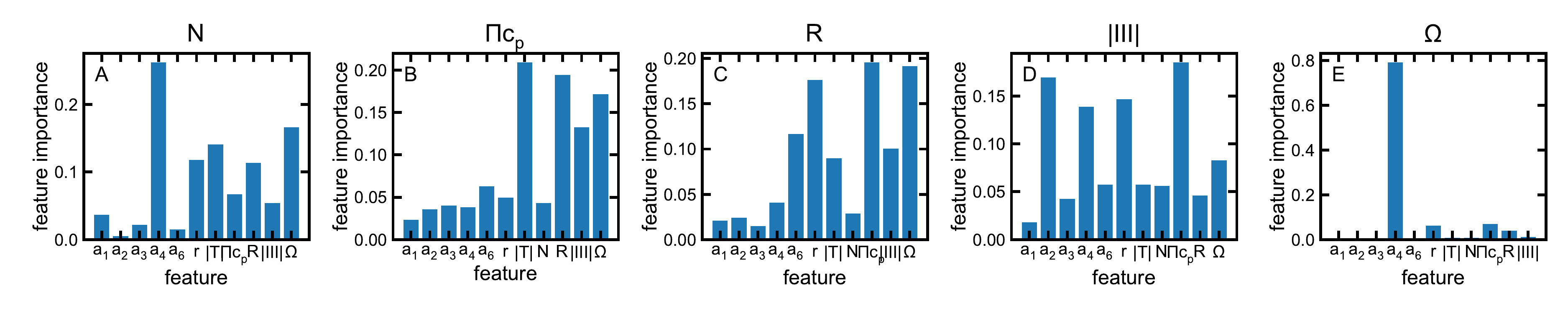}
\caption{\textbf{Feature importance considering as features all the quantities characterizing an elliptic curve.} Feature importance of the XGBoost regression models for predicting $N$ (Part A), $\prod\limits_{p \mid N} c_p$ (Part B), $R$ (Part C), $|\Sha|$ (Part D) and $\Omega$ (Part E).}
\label{regr_1}
\end{figure}

%%%%

We present the results in Table \ref{results_regression_3} of the accuracy measure NMAE by the 3 methods as in subsection \ref{ml1}: machine-learning by XGBoost, dummy regression and linear regression.
To read the table, of each of the 5 quantities given in a row, we use the other 4 to train the ML in order to predict this row.
We see that this is a significant improvement over the previous subsection and shows that, especially the XGBoost, the machine-learning can very confidently predict the Tamagawa number, the regulator and $|\Sha|$.
The feature importance is shown in Fig.~\ref{regr_1} and in table \ref{lin3}, we report the the statistics associated to the mixed predictions linear regression models. A comparison between the predictions of the linear models compared to XGboost is presented in Figure \ref{Figure_True_pred}.

\begin{table}[h!!!]
\begin{center}
\begin{tabular}{rrrrr}
\toprule
Predicted variable & R-squared & Adj. R-squared & F-statistic & Prob (F-statistic) \\ \hline
$N$ &      0.038 &           0.038 &         8999 &                   0 \\

$\prod\limits_{p \mid N} c_p$
 &      0.054 &           0.054 &        12960 &                   0 \\

$R$
 &      0.042 &           0.042 &         9889 &                   0 \\

$|\Sha|$
 &      0.017 &           0.017 &         3891 &                   0 \\

$\Omega$
 &      0.114 &           0.114 &        29110 &                   0 \\

\bottomrule
\end{tabular}
\end{center}
\caption{\textbf{Statistics of the linear regression models (mixed predictions).} We report the R squared, the adjusted R squared, the F statistics, and the p-value associated to the F statistics for the mixed predictions linear models. }
\label{lin3}
\end{table}

%%%
\begin{table}[h!!!]
    \centering
    \begin{tabularx}{\columnwidth}{@{}Xrrr@{}}
        \midrule
        \textbf{Quantity} & \textbf{ $F1_{micro}$ (XGBoost)} & \textbf{ $F1_{micro}$ (Dummy)} & \textbf{ $F1_{micro}$ (Logistic)} \\ \toprule

$r$ & $0.900\pm 0.001 $& $0.502\pm 0.001$& $0.730\pm 0.001$ \\ 
$|T|$ & $0.908\pm 0.001 $& $0.543\pm 0.001$& $0.567\pm 0.001$ \\
         \bottomrule
    \end{tabularx}
     \begin{tabularx}{\columnwidth}{@{}Xrrr@{}}
        \midrule
         & \textbf{ $F1_{macro}$ (XGBoost)} & \textbf{ $F1_{macro}$ (Dummy)} & \textbf{ $F1_{macro}$ (Logistic)} \\ \toprule

$r$ & $0.554\pm 0.001 $& $0.167\pm 0.001$& $0.387\pm 0.001$ \\ 
$|T|$ & $0.585\pm 0.015 $& $0.059\pm 0.001$& $0.090\pm 0.001$ \\ 
         \bottomrule
    \end{tabularx}
    
     \begin{tabularx}{\columnwidth}{@{}Xrrr@{}}
        \midrule
         & \textbf{ $MCC$ (XGBoost)} & \textbf{ $MCC$ (Dummy)} & \textbf{ $MCC$ (Logistic)} \\ \toprule
$r$ & $0.8311\pm 0.0005 $& $0.0000\pm 0.0000$& $0.5240\pm 0.0011$ \\ 
$|T|$ & $0.8302\pm 0.0014 $& $0.0000\pm 0.0000$& $0.1364\pm 0.0009$ \\ 
         \bottomrule
    \end{tabularx}
    \caption{\textbf{Performance of the classification models considering as features all the quantities characterizing an elliptic curve.} The scores $F1_{macro}$, $F1_{micro}$ and the Matthew correlation coefficient $MCC$, for XGboost (left column), the dummy regressor (central column) and a logistic regression (right column). The reported values are averages across $5$-fold cross-validations, with the corresponding standard deviations. Results are considerably improved compared to \cref{results_classification}.}
\label{results_classification3}
\end{table}
%%%

\begin{figure}[h!]
\centering
\includegraphics[width=\linewidth]{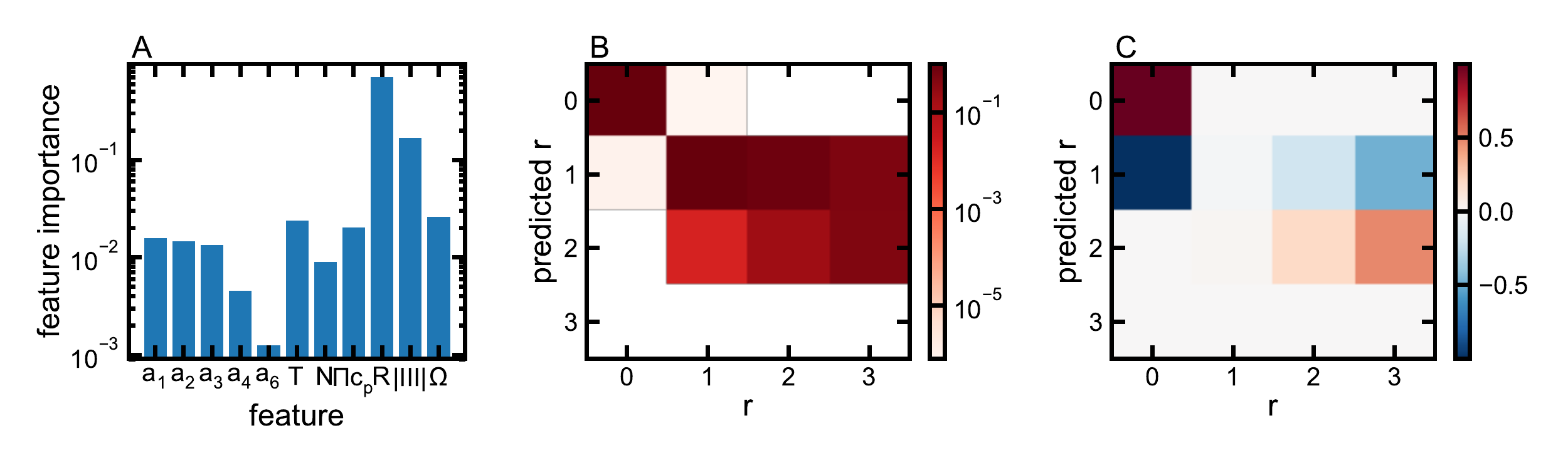}
\caption{\textbf{Prediction of $r$ considering as features all the quantities characterizing an ellipses.} (A) Importance of the different features to predict $r$. (B) Confusion matrix (normalized by column) showing the fraction of entries with rank $r$ with given $predicted$ $r$. (C) Difference between the confusion matrix obtained for the XGBoost and the dummy classifier. Results are averaged over a $5$-fold cross validation.}
\label{class_3}
\end{figure}

\begin{figure}[h!]
\centering
\includegraphics[width=\linewidth]{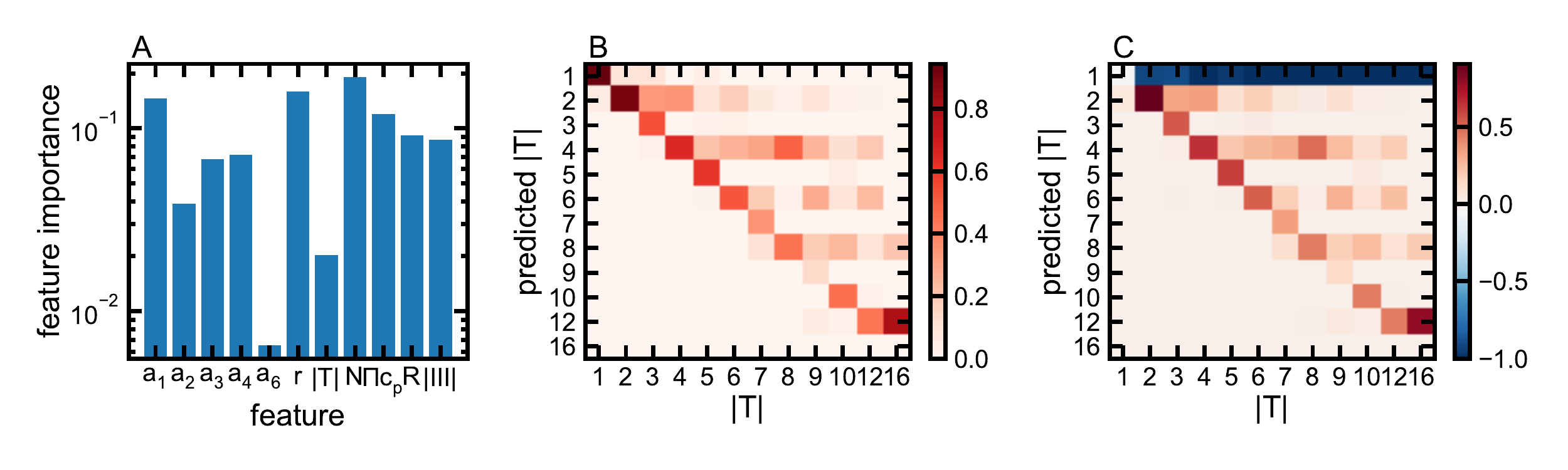}
\caption{\textbf{Prediction of $|T|$ considering as features all the quantities characterizing an elliptic curve.} (A) Importance of the different features to predict $|T|$. (B) Confusion matrix (normalized by column) showing the fraction of entries with value $|T|$ and given $predicted$ $|T|$. (C) Difference between the confusion matrix obtained for the XGBoost and the dummy classifier. Results are averaged over a $5$-fold cross validation.}
\label{class_4}
\end{figure}

Finally, let us use all quantities: the coefficients $a_i$ as well as $N$, $\prod\limits_{p \mid N} c_p$, $R$, $|\Sha|$, $\Omega$ to predict $r$ and $|T|$.
The accuracy measure $F1$ and $MCC$ (which should be close to 1 ideally)  are presented in Table \ref{results_classification3} and the feature importance, in Figures \ref{class_3} and \ref{class_4}.
We see that these are considerably improved compared to those obtained in  \cref{ml1} (cf. Tables \ref{results_regression_1} and \ref{results_classification}). This is somehow to be expected in light of the correlations observed in Figure \ref{correlation_matrix}.
In fact, even logistic regressions behave fairly well, with the $F1_{micro}$, $F1_{macro}$, and the $MCC$ scores subtiantially larger than those obtained by the dummy classifiers (see Table \ref{results_classification3}).
For reference, we include report the hyperplane equations of the linear regression models to see how each of the quantities can be fitted by all the others:
\begin{equation}
\begin{array}{rcl}
N &=& 195500.0000-1526.8435 a_1+     288.1786 a_2 -806.9174 a_3  -122.4328 a_4+      16.3690 a_6+   \\
	&&  8047.4199 r-10250.0000 |T|+    2192.0941 \prod\limits_{p \mid N} c_p+    3197.1580 R+    1293.3455 |\Sha|  -20330.0000 \Omega
\\
\prod\limits_{p \mid N} c_p&=&49.7747  +  14.9069 a_1+   3.6860 a_2+ 2.1424 a_3 -2.4121 a_4+   1.0613 a_6+ \\
&& 17.8877 r+  53.2012 |T|+   5.3471 N -14.5038 R  -4.3689 |\Sha|-27.7937 \Omega
\\
R&=&
  4.0689-0.0028 a_1 -0.1430 a_2 -0.2379 a_3 -0.2995 a_4+  0.1230 a_6+ \\
&&2.1850 r+  0.4585 |T|+  0.3910 N -0.7271 \prod\limits_{p \mid N} c_p -0.2167 |\Sha| -2.1082 \Omega
\\
|\Sha| &=&
  1.5946  +  0.0517 a_1+  0.0094 a_2 -0.0195 a_3 -0.6322 a_4+  0.4875 a_6 - \\
&&0.5472 r+  0.0112 |T|+  0.0756 N -0.1046 \prod\limits_{p \mid N} c_p -0.1035 R -0.3466 \Omega
\\
\Omega &=&
  0.6065 -0.0252 a_1+  0.0113 a_2+  0.0160 a_3 -0.0030 a_4+  0.0019 a_6+ \\
&& 0.1147 r -0.0717 |T| -0.0945 N -0.0530 \prod\limits_{p \mid N} c_p -0.0801 R -0.0276 |\Sha|
\end{array}
\end{equation}

%%

%%%%%%%%%%%%%%%%%%%%%%%%%
\section{Conclusions and Prospects}
In this paper, we initiated the study of the data science of the arithmetic of elliptic curves in light of the Birch-Swinnerton-Dyer Conjecture.
This is inspired by the the recent advances in the statistical investigation of Calabi-Yau manifolds, especially in the context of super-string theory \cite{He:2015fif,Altman:2018zlc}, as well as in the paradigm of machine-learning structures in geometry \cite{He:2017aed,He:2018jtw} and algebra \cite{He:2019nzx}.
Whilst we are still within the landscape of "Calabi-Yau-ness", it is expected that patterns in number theory should be much more subtle than those in geometry over $\IC$ and in combinatorics.
Nevertheless, BSD, residing at the crux between arithmetic and analysis, might be more susceptible to machine-learning and to pattern-recognition.

From our preliminary examinations on the extensive database of Cremona \cite{cremona,lmfdb} we have already found several interesting features.
First, we find that in the minimal Weierstra\ss\ representation, where a pair of coefficients $(a_4, a_6)$ clearly constitutes the principle component of the data, the distribution thereof follows a curious cross-like symmetry across rank, as shown in Figures \ref{f:sloga4a6} and \ref{a4-a6-pdf}. This is a highly-non-trivial symmetry since $a_{4,6} \leftrightarrow \pm a_{4,6}$ does not preserve rank.
This symmetry is reminiscent of mirror symmetry for Calabi-Yau threefolds.
In addition, the absence of data-points beyond the boundaries of the cross is also of note, much like that Hodge plot for the Calabi-Yau threefolds.

Over all, the distribution of the Euclidean distance of $(a_4, a_6)$ to the origin, as well as that of the RHS of the Strong BSD, viz., the quantity $\frac{|\Sha| \cdot \Omega \cdot R \cdot \prod\limits_{p \mid N} c_p}{|T|^2}$ (cf.~conjectures in Sec.~\ref{s:conj}), are best described by a Beta-distribution, which is selected in both cases among 85 continuous distributions using the Akaike Information Criterion.
Organized by rank, these distributions also vary.

One further visualize the data, the tuples consisting of the coefficients $(a_1,a_2,a_3,a_4,a_6)$ as well as the BSD tuple $(N,\ |T|,\ \prod\limits_{p \mid N} c_p,\ \Omega,\ R, \Sha)$ for ranks $r =0,1,2$ using the standard techniques from topological data analysis.
The bar-codes are shown in Figures \ref{f:coefbarcode} and \ref{f:6barcode}.
While the Weierstra\ss coefficients show little variation over rank, the BSD tuple does show differences over $r$.
Moreover, as expected, the divisibility of the conductor $N$ influences the barcodes.

Finally, emboldened by the recent success in using machine-learning to computing bundle cohomology on algebraic varieties without recourse to sequence-chasing and Gr\''obner bases as well as recognizing whether a finite group is simple directly by looking at the Cayley table, we asked the question of whether one can ``predict'' quantities otherwise difficult to compute directly from ``looking'' at the shape of the elliptic curve.
Ideally, one would have hoped that training on $a_i$, one could predict any of the BSD quantities to high precision, as was in the cohomology case.
However, due to the very high variation in the size of $a_i$, one could not find a good machine-learning technique, decision trees, support-vector machines or neural networks, which seems to achieve this.
This is rather analogous to the (expected) failure of predicting prime numbers using AI.
Nevertheless, the BSD quantities, when {\it mixed} with the Weierstra\ss\ coefficients, does behave well under machine-learning.
For instance, the Matthew correlation coefficient between predicted and true values of $r$ and $|T|$ is $\sim 0.83$.

At some level, the experiments here, in conjunction with those in \cite{He:2017aed,Krefl:2017yox,Ruehle:2017mzq,Carifio:2017bov,Bull:2018uow,Jejjala:2019kio,He:2019nzx,Brodie:2019dfx,Ashmore:2019wzb}, tend to show a certain {\it hierarchy of difficulty} in how machine-learning responds to problems in mathematics.
Understandably, number theory is the most difficult: as a reprobate, \cite{He:2017aed} checked that trying to predict the next prime number, for instance, seems unfeasible for simple neural networks.
On the other hand, algebraic geometry over the complex numbers seems to present a host of amenable questions such as bundle cohomology, recognition of elliptic fibrations or calculating Betti numbers.
In between lie algebra and combinatorics, such as the structure of finite groups, where precision/confidence of the cross-validation is somewhat intermediate.
It is therefore curious that in this paper one sees that a problem such as BSD, which resides in between arithmetic and geometry, is better behaved under machine-learning than a direct attack on patterns in primes.

%%%%%%%%%%%%%%%%%%%%%
%%%=============================================================
\section*{Acknowledgments}
YHH would like to thank the Science and Technology Facilities Council,
UK, for grant ST/J00037X/1, Nankai University, China for a chair professorship and Merton College, Oxford, for enduring support.

\pagebreak

\begin{appendices}
%\section{XGboost hyperparameters optimization.}
%\label{par_opt}
%To prevent overfitting, control the complexity of the model and make it robust to noise, we tune the following parameters. The other parameters are kept to their default value as per the XGboost classifier implementation in Python \cite{chen2016xgboost}.
%
%\begin{itemize}
%\item{\textbf{$n\_estimators$: }}{The number of trees (between $150$ and $500$).}
%\item{\textbf{$learning\_rate$: }}{It controls how the output of each tree contributes to update the estimates (between $0.01$ and $0.07$).}
%\item{\textbf{$subsample$: }}{The fraction of observations to be selected for each tree (included between $0.3$ and $1$). }
%\item{\textbf{$max\_depth$: }}{The maximum depth of a tree (included between $3$ and $9$).}
%\item{\textbf{$min\_child\_weight$: }}{Defines the minimum sum of weights of all observations required in a child (included between $1$ and $10$).}
%\item{\textbf{$gamma$: }}{A node is split only when the resulting split gives a positive reduction in the loss function. Gamma specifies the minimum loss reduction required to make a split (included between $10.5$ and $5$).}
%\item{\textbf{$colsample\_bytree$: }}{Denotes the fraction of columns to be randomly samples for each tree (included between $0.8$ and $1$).}
%\end{itemize}
%
%Parameters are optimized against the average F1-score across ranks, using a 10-fold cross validation.
%
%We obtain the following values: $n\_estimators=420$, $colsample\_bytree=0.8$, $gamma=1$, $learning\_rate=0.0285$, $max\_depth=9$,$ min\_child\_weight=5$, $subsample=1$.

\section{Learning Curves}
\label{learning_curve}
\begin{figure}[h!]
\centering
\includegraphics[width=\linewidth]{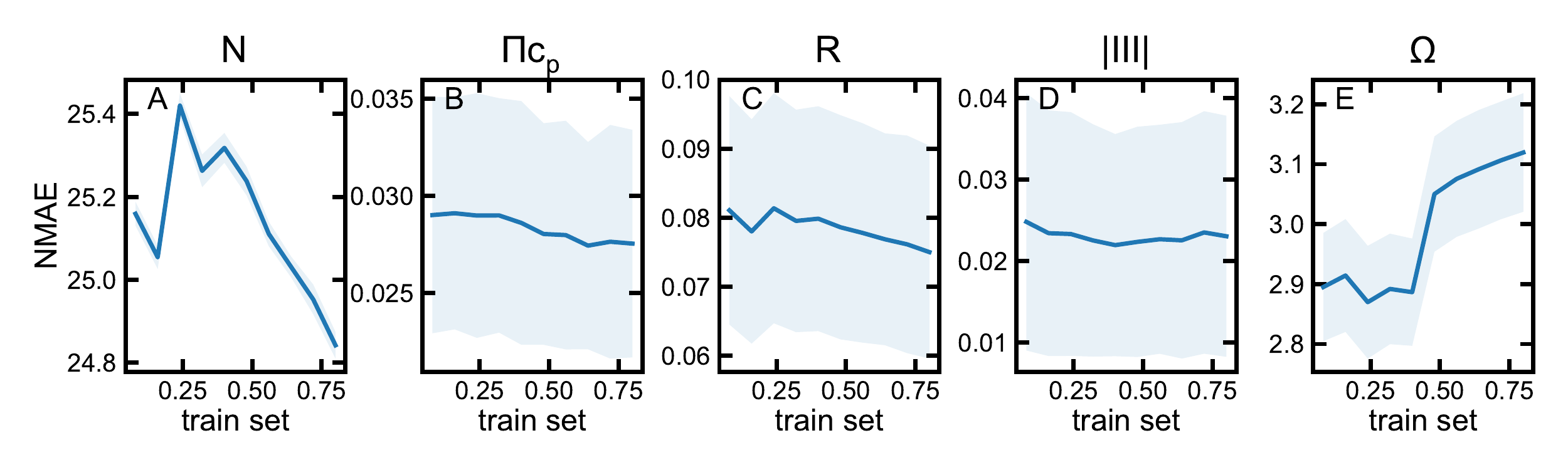}
\caption{\textbf{Learning curves for the regression models.} The Normalized Median Absolute Error as a function of the train set size of the XGBoost regression models for predicting $N$ (A), $\prod\limits_{p \mid N} c_p$ (B), $R$ (C), $|\Sha|$ (D) and $\Omega$ (E). The shaded areas correspond to standard deviation across a $5-fold$ cross validation. }
\label{learning_curves_1}
\end{figure}

\begin{figure}[h!]
\centering
\includegraphics[width=\linewidth]{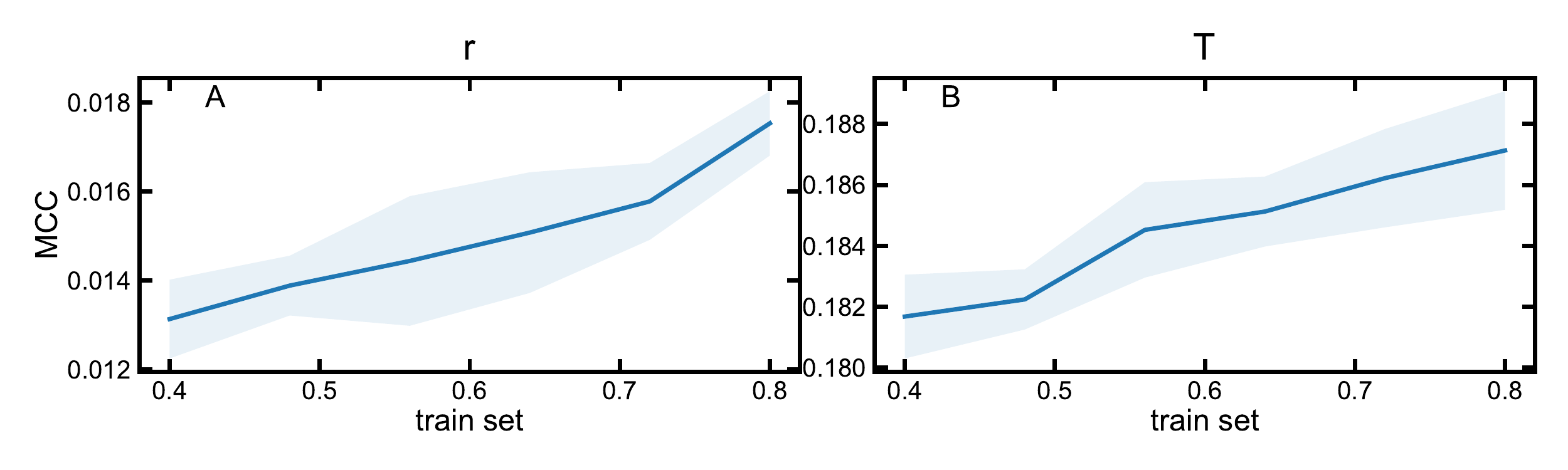}
\caption{\textbf{Learning curves for the classification models.} The Matthew coefficient as a function of the train set size of the XGBoost regression models for predicting $r$ (A) and $|T|$ (B). The shaded areas correspond to standard deviation across a $5-fold$ cross validation. }
\label{learning_curves_2}
\end{figure}

To prevent overfitting, we compute the learning curves of the regression (figure \ref{learning_curves_1}) and classification (figure \ref{learning_curves_2}) models in \cref{ml1}. We find that $80\%$ of the data is a good choice for the training set size, suggesting a 5-fold cross validation. 

%%%%%%%%%%%%%%%
\section{Comparison with SVM}\label{svm}
In this section, we compare the performance of the  XGBoost models with Support Vector Machine (SVM) models. SVM models are very long to train, hence we focus, for this task, on a subset of $100,000$ examples. In table \ref{svm1} we report the performance of the regression models used to predict 
 $N$, $\prod\limits_{p \mid N} c_p$, $R$, $|\Sha|$ and $\Omega$. Only in the case of $\Omega$, the SVM model performs better than XGBoost.

\begin{table}[ht]
    \centering
    \begin{tabularx}{\columnwidth}{@{}Xrrr@{}}
        \midrule
        \textbf{Quantity} & \textbf{NMAE (XGBoost)} & \textbf{NMAE (Dummy)} & \textbf{NMAE (SVM)} \\ \toprule
$N$ & $114881.964\pm 364.409 $& $115664.993\pm 384.936$& $115690.292\pm 417.181$ \\ 
$\prod\limits_{p \mid N} c_p$ & $278.372\pm 34.807 $& $273.775\pm 27.528$& $275.182\pm 27.412$ \\ 
$R$ & $17.493\pm 4.178 $& $15.124\pm 4.137$& $15.417\pm 4.067$ \\ 
$|\Sha|$ & $4.938\pm 1.156 $& $4.868\pm 1.223$& $4.893\pm 1.218$ \\ 
$\Omega$ & $0.498\pm 0.004 $& $0.584\pm 0.005$& $0.607\pm 0.005$ \\ 
         \bottomrule
    \end{tabularx}
    \caption{\textbf{Performance of the regression models.} The normalized median absolute error $NMAE$, for XGboost (left column), the dummy regressor (central column) and Support Vector Machine Regression (right column). The reported values are averages across $5$-fold cross-validations, with the corresponding standard deviations.}
\label{svm1}
\end{table}

%%%%%%%%%%%%%%%%%%%%%%%%
\section{Characteristics of the Weierstra\ss\  coefficients.}
\renewcommand{\arraystretch}{0.7}

\begin{footnotesize}
\begin{longtable}{rrrrrrrrr}
\toprule
a1 & a2 & a3 & rank &      size &     $\overline{a_4}$    &     $s_{a_4}$   &    median     &            zero entries \\
\midrule
 0 & -1 &  0 &    0 &  126135 & -5E+09 & 7E+11 & -8E+03 &          98 \\
 1 &  1 &  1 &    0 &   67834 & -8E+10 & 7E+12 & -2E+04 &          54 \\
 1 &  1 &  0 &    0 &   69759 & -2E+11 & 3E+13 & -1E+04 &          47 \\
 1 &  0 &  1 &    0 &   71309 & -9E+10 & 1E+13 & -2E+04 &          35 \\
 1 &  0 &  0 &    0 &   66411 & -1E+12 & 2E+14 & -2E+04 &          42 \\
 1 & -1 &  1 &    0 &   96995 & -1E+11 & 2E+13 & -2E+04 &          41 \\
 0 &  1 &  1 &    0 &   18016 & -4E+10 & 3E+12 & -2E+03 &          38 \\
 0 &  1 &  0 &    0 &  118942 & -1E+10 & 1E+12 & -9E+03 &         108 \\
 0 &  0 &  1 &    0 &   28440 & -1E+11 & 2E+13 & -3E+03 &         546 \\
 1 & -1 &  0 &    0 &  102769 & -2E+11 & 5E+13 & -2E+04 &          97 \\
 0 &  0 &  0 &    0 &  172238 & -6E+09 & 9E+11 & -1E+04 &         832 \\
 0 & -1 &  1 &    0 &   17238 & -1E+10 & 1E+12 & -2E+03 &          40 \\
 0 & -1 &  1 &    1 &   24593 & -1E+09 & 1E+11 & -1E+03 &          65 \\
 1 & -1 &  0 &    1 &  127198 & -1E+11 & 2E+13 & -1E+04 &         150 \\
 1 &  0 &  0 &    1 &   98092 & -5E+11 & 1E+14 & -7E+03 &          77 \\
 0 &  1 &  1 &    1 &   27360 & -5E+10 & 4E+12 & -1E+03 &          54 \\
 1 &  0 &  1 &    1 &   94595 & -3E+11 & 6E+13 & -7E+03 &          62 \\
 1 & -1 &  1 &    1 &  128957 & -2E+10 & 2E+12 & -9E+03 &          40 \\
 0 &  0 &  0 &    1 &  213780 & -5E+10 & 2E+13 & -6E+03 &         962 \\
 0 &  1 &  0 &    1 &  157003 & -8E+09 & 1E+12 & -5E+03 &         164 \\
 1 &  1 &  0 &    1 &   88403 & -5E+10 & 4E+12 & -7E+03 &         107 \\
 0 & -1 &  0 &    1 &  155604 & -1E+10 & 2E+12 & -5E+03 &         159 \\
 0 &  0 &  1 &    1 &   39235 & -5E+10 & 7E+12 & -2E+03 &         608 \\
 1 &  1 &  1 &    1 &   91717 & -3E+10 & 4E+12 & -7E+03 &          87 \\
 1 &  1 &  0 &    2 &   18293 & -2E+08 & 2E+10 & -1E+03 &          28 \\
 1 &  0 &  0 &    2 &   25286 & -5E+07 & 2E+09 & -2E+03 &          23 \\
 1 & -1 &  1 &    2 &   28940 & -5E+07 & 6E+09 & -2E+03 &          17 \\
 0 & -1 &  0 &    2 &   30236 & -3E+07 & 2E+09 & -1E+03 &          44 \\
 1 &  0 &  1 &    2 &   20907 & -2E+08 & 2E+10 & -1E+03 &          17 \\
 0 &  0 &  0 &    2 &   40731 & -6E+07 & 6E+09 & -2E+03 &         126 \\
 1 & -1 &  0 &    2 &   25793 & -6E+08 & 8E+10 & -2E+03 &          46 \\
 1 &  1 &  1 &    2 &   21197 & -7E+07 & 5E+09 & -1E+03 &          23 \\
 0 &  0 &  1 &    2 &   11187 & -2E+07 & 9E+08 & -5E+02 &          96 \\
 0 &  1 &  1 &    2 &    9609 & -1E+08 & 1E+10 & -5E+02 &          19 \\
 0 & -1 &  1 &    2 &    7582 & -7E+06 & 2E+08 & -3E+02 &          22 \\
 0 &  1 &  0 &    2 &   34585 & -2E+07 & 9E+08 & -1E+03 &          36 \\
 1 & -1 &  0 &    3 &     551 & -8E+04 & 1E+06 & -3E+02 &           1 \\
 0 &  0 &  0 &    3 &     698 & -2E+04 & 2E+05 & -3E+02 &           0 \\
 1 &  1 &  0 &    3 &     496 & -2E+04 & 2E+05 & -2E+02 &           0 \\
 0 &  0 &  1 &    3 &     506 & -2E+04 & 2E+05 & -2E+02 &           2 \\
 0 & -1 &  1 &    3 &     399 & -8E+04 & 1E+06 & -2E+02 &           1 \\
 0 &  1 &  0 &    3 &     722 & -8E+03 & 4E+04 & -4E+02 &           1 \\
 0 & -1 &  0 &    3 &     659 & -1E+04 & 6E+04 & -3E+02 &           0 \\
 1 &  0 &  0 &    3 &     612 & -6E+03 & 3E+05 & -4E+02 &           1 \\
 0 &  1 &  1 &    3 &     426 &  4E+03 & 9E+05 & -3E+02 &           1 \\
 1 & -1 &  1 &    3 &     604 & -1E+04 & 7E+04 & -4E+02 &           3 \\
 1 &  0 &  1 &    3 &     548 & -1E+04 & 1E+05 & -3E+02 &           3 \\
 1 &  1 &  1 &    3 &     458 & -2E+04 & 4E+05 & -2E+02 &           0 \\
 1 & -1 &  0 &    4 &       1 & -8E+01 &   NAN & -8E+01 &           0 \\
\bottomrule
\caption{For given values of $a_1$,$a_2$,$a_3$, and $r$, the table reports the number of elliptic curves (size), and some statistics of the Weierstra\ss\  coefficient $a_4$ including the mean ($\overline{a_4}$), the standard deviation ($s_{a_4}$), the median (median) and the number of zero entries (zero entries)}
\label{Table1}
\end{longtable}
\end{footnotesize}

%%%%%%%%%%%%%
\begin{footnotesize}
\begin{longtable}{rrrrrrrrr}
\toprule
a1 & a2 & a3 & rank &      size &     $\overline{a_6}$    &     $s_{a_6}$   &    median     &            zero entries \\
\midrule
 0 & -1 &  0 &    0 &  126135 & -9E+15 & 4E+18 & -5E+02 &         217 \\
 1 &  1 &  1 &    0 &   67834 & -3E+17 & 9E+19 & -1E+03 &           1 \\
 1 &  1 &  0 &    0 &   69759 & -3E+18 & 6E+20 & -7E+02 &         202 \\
 1 &  0 &  1 &    0 &   71309 &  2E+17 & 2E+20 & -2E+03 &           2 \\
 1 &  0 &  0 &    0 &   66411 & -4E+19 & 1E+22 & -5E+03 &         176 \\
 1 & -1 &  1 &    0 &   96995 & -6E+17 & 3E+20 & -3E+03 &           1 \\
 0 &  1 &  1 &    0 &   18016 & -2E+17 & 3E+19 & -4E+02 &           2 \\
 0 &  1 &  0 &    0 &  118942 & -2E+16 & 6E+18 & -1E+03 &         226 \\
 0 &  0 &  1 &    0 &   28440 & -3E+18 & 5E+20 & -6E+02 &           1 \\
 1 & -1 &  0 &    0 &  102769 &  7E+18 & 2E+21 & -8E+02 &         206 \\
 0 &  0 &  0 &    0 &  172238 &  1E+16 & 6E+18 & -8E+02 &         486 \\
 0 & -1 &  1 &    0 &   17238 & -1E+16 & 4E+18 & -2E+02 &           1 \\
 0 & -1 &  1 &    1 &   24593 &  1E+15 & 2E+17 &  1E+01 &          17 \\
 1 & -1 &  0 &    1 &  127198 & -2E+18 & 6E+20 & -4E+00 &         271 \\
 1 &  0 &  0 &    1 &   98092 & -3E+19 & 9E+21 &  5E+01 &         246 \\
 0 &  1 &  1 &    1 &   27360 & -1E+17 & 4E+19 & -2E+00 &          18 \\
 1 &  0 &  1 &    1 &   94595 &  1E+19 & 3E+21 &  2E+01 &          27 \\
 1 & -1 &  1 &    1 &  128957 & -2E+14 & 1E+19 &  4E+01 &          15 \\
 0 &  0 &  0 &    1 &  213780 & -1E+18 & 6E+20 & -5E-01 &         604 \\
 0 &  1 &  0 &    1 &  157003 & -3E+15 & 5E+18 &  2E+01 &         281 \\
 1 &  1 &  0 &    1 &   88403 & -1E+17 & 4E+19 &  0E+00 &         271 \\
 0 & -1 &  0 &    1 &  155604 &  2E+15 & 1E+19 &  0E+00 &         301 \\
 0 &  0 &  1 &    1 &   39235 &  6E+17 & 1E+20 & -2E+01 &          24 \\
 1 &  1 &  1 &    1 &   91717 &  1E+17 & 4E+19 &  7E+00 &          17 \\
 1 &  1 &  0 &    2 &   18293 &  8E+13 & 1E+16 &  1E+03 &          42 \\
 1 &  0 &  0 &    2 &   25286 &  4E+12 & 4E+14 &  1E+04 &          42 \\
 1 & -1 &  1 &    2 &   28940 & -1E+13 & 3E+15 &  1E+04 &          21 \\
 0 & -1 &  0 &    2 &   30236 &  7E+11 & 3E+14 &  2E+03 &          58 \\
 1 &  0 &  1 &    2 &   20907 &  7E+13 & 9E+15 &  3E+03 &          22 \\
 0 &  0 &  0 &    2 &   40731 & -2E+13 & 3E+15 &  4E+03 &          82 \\
 1 & -1 &  0 &    2 &   25793 &  7E+14 & 1E+17 &  2E+03 &          57 \\
 1 &  1 &  1 &    2 &   21197 & -9E+12 & 1E+15 &  5E+03 &          23 \\
 0 &  0 &  1 &    2 &   11187 &  1E+12 & 1E+14 &  1E+03 &          26 \\
 0 &  1 &  1 &    2 &    9609 & -4E+13 & 4E+15 &  2E+03 &          18 \\
 0 & -1 &  1 &    2 &    7582 &  4E+10 & 8E+12 &  6E+02 &          24 \\
 0 &  1 &  0 &    2 &   34585 &  5E+11 & 9E+13 &  5E+03 &          42 \\
 1 & -1 &  0 &    3 &     551 &  1E+08 & 3E+09 &  2E+03 &           0 \\
 0 &  0 &  0 &    3 &     698 & -2E+06 & 2E+08 &  2E+03 &           0 \\
 1 &  1 &  0 &    3 &     496 &  3E+06 & 9E+07 &  9E+02 &           1 \\
 0 &  0 &  1 &    3 &     506 & -6E+06 & 2E+08 &  1E+03 &           4 \\
 0 & -1 &  1 &    3 &     399 &  1E+08 & 3E+09 &  6E+02 &           3 \\
 0 &  1 &  0 &    3 &     722 &  1E+06 & 1E+07 &  3E+03 &           0 \\
 0 & -1 &  0 &    3 &     659 &  3E+06 & 3E+07 &  2E+03 &           0 \\
 1 &  0 &  0 &    3 &     612 & -5E+06 & 2E+08 &  3E+03 &           1 \\
 0 &  1 &  1 &    3 &     426 &  1E+08 & 2E+09 &  1E+03 &           4 \\
 1 & -1 &  1 &    3 &     604 &  3E+06 & 4E+07 &  3E+03 &           4 \\
 1 &  0 &  1 &    3 &     548 &  5E+06 & 6E+07 &  2E+03 &           4 \\
 1 &  1 &  1 &    3 &     458 &  5E+07 & 1E+09 &  1E+03 &           2 \\
 1 & -1 &  0 &    4 &       1 &  3E+02 &   NAN &  3E+02 &           0 \\
\bottomrule

\caption{For given values of $a_1$,$a_2$,$a_3$, and $r$, the table reports the number of curves (size), and some statistics of the Weierstra\ss\  coefficient $a_6$ including the mean ($\overline{a_6}$), the standard deviation ($s_{a_6}$), the median (median) and the number of zero entries (zero entries)}
\label{Table2}
\end{longtable}
\end{footnotesize}

\end{appendices}

{\small
%%
%%%%%%%%%%%%%%%%%%%%%%%%%%%%=======================================

}

\end{document}